\documentclass[final]{SIAM}

\usepackage{graphicx}%
\usepackage{multirow}%
\usepackage{amsmath,amssymb,amsfonts}%
\usepackage{xcolor}%
\usepackage{textcomp}%
\usepackage{manyfoot}%
\usepackage{booktabs}%
\usepackage{algorithm}%
\usepackage{algorithmicx}%
\usepackage{algpseudocode}%
\usepackage{listings}%

\usepackage{mathptmx}

\newsiamthm{remark}{Remark}

\usepackage{array}
\newcolumntype{L}{>{\centering\arraybackslash}m{3cm}}

\def\<{{\langle}} 
\def\>{{\rangle}} 

\newcommand{\vertiii}[1]{{\left\vert\kern-0.25ex\left\vert\kern-0.25ex\left\vert #1 
    \right\vert\kern-0.25ex\right\vert\kern-0.25ex\right\vert}}

\title{Strong convergence with error estimates for a stochastic compartmental model of electrophysiology 
}

\date{\today}

\usepackage[sc]{mathpazo}
\usepackage[T1]{fontenc}
\usepackage[mmddyyyy]{datetime}
\usepackage{advdate}
\newdateformat{syldate}{\twodigit{\THEMONTH}/\twodigit{\THEDAY}}
\newsavebox{\MONDAY}\savebox{\MONDAY}{Mon}
\newcommand{\week}[1]{%
  \paragraph*{\kern-2ex\quad #1, \syldate{\today} - \AdvanceDate[4]\syldate{\today}:}
  \ifdim\wd1=\wd\MONDAY
    \AdvanceDate[7]
  \else
    \AdvanceDate[7]
  \fi%
}
\usepackage{setspace}
\usepackage{multicol}
\usepackage{url}
\usepackage{hyperref}
\usepackage{lastpage}
\usepackage{amsmath}
\usepackage{layout}

\usepackage{array, xcolor}
\usepackage{color,hyperref}
\definecolor{clemsonorange}{HTML}{EA6A20}
\hypersetup{colorlinks,breaklinks,linkcolor=clemsonorange,urlcolor=clemsonorange,anchorcolor=clemsonorange,citecolor=black}

\title{Isochronal Phase Reduction and Speed Correction of a Pulse in a Stochastic Kinematic Model 
}

\begin{document}

\title{Isochronal Phase Reduction and Speed Correction of a Pulse in a Stochastic Kinematic Model 
}

\author{
Joshua A. McGinnis\footnotemark[2] \thanks{Corresponding Author (\email{jam887@sas.upenn.edu}})
\and
Xinbo Li\thanks{Department of Mathematics, University of Pennsylvania, Philadelphia, PA, 19104, USA.}
 \and Toshiyuki Ogawa\thanks{School of Interdisciplinary Mathematical Sciences, Meiji University, 4-21-1 Nakano, Nakanoku, Tokyo 164-8525, Japan.}
 \and
 Yoichiro Mori \footnotemark[2] \thanks{Department of Biology, University of Pennsylvania, Philadelphia, PA, 19104, USA. }
}

\maketitle


\begin{abstract}
We develop a method for computing the stochastic wave speed of pulse solutions in kinematic equations subject to small stochastic forcing based on the isochronal phase reduction. These kinematic equations arise as the singular limit of sharp pulse solutions in the FitzHugh–Nagumo system, and our approach contributes a new perspective and method to the growing body of work on stochastic wave propagation in excitable media. The method yields an effective Itô process for the wave’s position. The coefficients of the Itô process can be computed deterministically allowing for efficient computation. We demonstrate the efficiency and accuracy of our method through numerical demonstrations. 
\end{abstract}


\section{Introduction}
In this paper, we present a numerical procedure for computing the stochastic wave speed of a traveling pulse solution in kinematic equations subject to small noise. These kinematic equations arise as a singular limit of the FitzHugh–Nagumo equations, capturing the essential dynamics of propagation for sharp pulses. The procedure involves computing the coefficients of an Itô process that approximates the \textbf{isochronal phase}, which we adopt as our definition of the stochastic phase. Importantly, the procedure is fully deterministic: it requires solving a fixed number of PDEs, with the achievable accuracy increasing as the noise amplitude $\sigma$ decreases. This allows us to compute the mean and variance of the stochastic phase directly from the coefficients, bypassing the need for computationally intensive Monte Carlo simulations of the underlying SPDE. As a result, our approach offers a significantly more efficient alternative for quantifying the stochastic wave speed.

We begin with some background on the derivation of the kinematic equations. Our starting point is the Fitzhugh-Nagumo system,
\begin{equation}
\begin{aligned}
\epsilon\partial_t u(t,x) &=  \epsilon^2\partial_x^2u(t,x) + f_\text{cub}(u(t,x)) -w(t,x)
\\
\partial_t w(t,x) &= \delta \partial_x^2 w(t,x) +\beta u(t,x) -\gamma w(t,x),
\end{aligned}
\label{PDE:Fitz}
\end{equation}
posed on a periodic domain, $x \in \mathbb{R}/(L\mathbb{Z})$, where $f_\text{cub}$ is a cubic non-linearity with bistability. Often $f$ has the form $f(u) = u(1-u)(u-\alpha)$ with $\alpha$ the so called tuning parameter. This system is known to have stable, time periodic solutions consisting of one or multiple pulses moving with some wave speed. In this paper, we focus on the case of a single pulse denoted by $(u_{\epsilon,*}(x-c_\epsilon t),w_{\epsilon,*}(x-c_\epsilon t))$ moving with speed $c_\epsilon$. To fix ideas, we always consider the case where $c_\epsilon>0$ and a pulse moving to the right.

As $\epsilon \to 0$, the pulse becomes surrounded by narrow boundary layers and $c_\epsilon$ approaches a limit $c_0$. As a result the pulse sharpens. More specifically, in this singular limit, the spatial domain divides into two contiguous regions $E(t) \subset \mathbb{R}/(L\mathbb{Z})$ and $R(t)\subset \mathbb{R}/(L\mathbb{Z})$ representing the excited and the relaxation portions of the pulse solution, respectively.

These regions are separated by narrow boundary layers at $x_-(t)$ and $x_+(t)$, which, for an unperturbed pulse, move at constant speed $c_\epsilon$: $x_\pm(t) = c_\epsilon t+x_{\pm}(0)$. To understand the structure of the pulse near these interfaces, we consider the system in a moving frame, so that $\partial_t u_{\epsilon,*} = -c_\epsilon\partial_x u_{\epsilon,*}$, and the first equation becomes
\[ -\epsilon c_\epsilon \partial_x u_{\epsilon,*}   = \epsilon^2 \partial_{xx}^2u_{\epsilon,*}+f_{\text{cub}}(u_{\epsilon,*})-w_{\epsilon,*}.\]
In regions where the solution varies slowly in space, both derivative terms on the left- and right-hand sides are small compared to the reaction terms, and the leading-order balance becomes algebraic: $f_{\text{cub}}(u_{\epsilon,*}) \approx w_{\epsilon,*}$. Since $f_{\text{cub}}$ is bistable, this relation defines two stable branches of $u$, which we  denote by $\widetilde{g}_E(w_{\epsilon, *})$ and $\widetilde{g}_R(w_{\epsilon,*})$, corresponding to the excited and relaxation regions.

The pulse solution switches between these two branches across the boundary layers, where the spatial derivatives become large enough that the drift and diffusion terms can no longer be neglected. Hence,
 \begin{equation}
   \begin{cases} 
   u_{\epsilon,*} \approx \widetilde{g}_E(w_{\epsilon,*}) & x \in E(t) \\
   u_{\epsilon,*} \approx \widetilde{g}_R(w_{\epsilon,*}) & x \in R(t).
   \end{cases}
 \end{equation}
The approximation becomes equality in the singular limit $\epsilon \to 0$. In that limit, the profile  $w_*=w_{0,*}$  satisfies a piecewise equation
\begin{equation}
\partial_t w_* = \delta\partial_x^2 w_* +\begin{cases}
g_E(w_*) & x \in E(t) \\
g_R(w_*) & x \in R(t),
\end{cases}
\end{equation}
where $g_E$ and $g_R$ result from substituting the respective branches into the second equation in \eqref{PDE:Fitz}.
 Moreover, the wave speed satisfies a matching condition across the interface, encoded by a function $c(w)$  s.t. $\pm c(w_*(x_\pm-c_0t)) =c_0.$ Details regarding this function can be found in \cite{Keener1}. This motivates the study of the \textbf{kinematic} equations 
\begin{equation*}
\begin{aligned}
    \partial_t w &= \delta \partial_x^2 w + \begin{cases}
g_E(w) & x \in E(x_-,x_+) \\
g_R(w) & x \in R(x_-,x_+)
\end{cases}
\\
\dfrac{dx_\pm}{dt} &=\pm c(w(x\pm)),
\end{aligned}
\end{equation*}
which describe the approximate behavior of sharp traveling pulse solution in the FitzHugh-Nagumo systems when $\epsilon\ll1.$ 

We are interested in cases of multiplicative stochastic forcing but of a particular form. Our starting point is the SPDE
\begin{equation}
\begin{aligned}
\epsilon d u(t,x) &=  \left(\epsilon^2\partial_x^2u(t,x) + f_\text{cub}(u(t,x)) -w(t,x)\right)dt
\\
d w(t,x) &= \left(\delta \partial_x^2 w(t,x) +\beta u(t,x) -\gamma w(t,x)\right)dt + \sigma \widetilde{h}(w ,u)dW(t,x),
\end{aligned}
\label{SPDE:Fitz}
\end{equation}
 where $dW(t,x)$ denotes translationally invariant, spatially correlated noise. 
 
When $\epsilon>0, \delta>0$, and $\sigma$ is small there is well-developed theory for calculating a stochastic wave speed for a pulse solution of \eqref{SPDE:Fitz}. In particular, the work of Hamster and Hupkes 
\cite{Hamster_Hupkes_1,hamster_hupkes_2,Hamster_Hupkes_3} rigorously develops a method for extracting high order approximations of a stochastic phase from pulse solutions using the spectral and stability theory of the system. Their framework includes stochastic forcing in both equations, including the equation for $u$.

 In contrast, we focus on a different parametric regime. In many applications, the parameters satisfy $\epsilon \ll\sigma^2$, while $\delta$ may be comparable to or even smaller than $\sigma^2$, and in some cases, $\delta = 0$. Moreover, we consider forcing only in the equation for $w$, not for $u$. This setup arises naturally in models where $w$ is derived from an underlying microscopic stochastic process and is the main source of intrinsic randomness.

For example, in models of neuronal action potentials, $u$ represents membrane voltage, and $w$ represents ion channel activity. Since ion channels are typically not directly coupled across space, it is natural to set $\delta = 0$. In such systems, the primary source of noise arises from the stochastic opening and closing of ion channels; see, e.g., \cite{Austin}. Whether $\epsilon$ is in fact small in such biological systems is debatable, but we nevertheless believe our parameter regime is both mathematically natural and physically relevant.

The singular limit $\epsilon \to 0$ of stochastic pulse solutions arising in \eqref{SPDE:Fitz} has not yet been rigorously justified. Nonetheless, one formally expects that this limit yields a stochastically forced version of the kinematic equations:
\begin{equation*}
\begin{aligned}
   d w &= \delta \partial_x^2 wdt + \begin{cases}
g_E(w)dt+\sigma h_E(w)dW(t,x) & x \in E(x_-,x_+) \\
g_R(w)dt + \sigma h_R(w)dW(t,x) & x \in R(x_-,x_+)
\end{cases}
\\
dx_{\pm} &=\pm c(x_{\pm}(t))dt,
\end{aligned}
\end{equation*}
which will be the equations we work with throughout the remainder of the paper. In particular, we define a procedure for calculating the stochastic wave speed in such systems. Whether $\delta$ is positive or $0$ makes little difference for our method, and since the expressions we obtain in the sequel are quite cumbersome, we set $\delta=0$.

Now we more precisely describe the isochronal map and describe our main contribution. Consider a deterministic system with a stable, one-dimensional center manifold. Given initial conditions in a neighborhood of this manifold, one may ask: where does the trajectory converge on the manifold as $t \to \infty$? The mapping that assigns each such nearby point to its asymptotic limit on the center manifold is known as the \textbf{isochronal map}.

In the presence of stochastic forcing, the system does not settle on the manifold, but instead fluctuates in its vicinity. For sufficiently small noise, the transverse displacement from the manifold behaves approximately like an Ornstein–Uhlenbeck process. As a result, trajectories are expected to remain close to the manifold and within its region of stability for exponentially long times.

This allows us to define what we call the \textbf{stochastic phase}: at any given time $t_0$, it is the image under the isochronal map of the system’s current state. Equivalently, the stochastic phase at time $t_0$ is the location on the center manifold where the trajectory would converge if the noise were turned off at that moment. In practice, this phase is typically identified with a coordinate which may be used to parametrize the center manifold.

While the stochastic phase is well-defined via the isochronal map, computing it directly is generally impractical, as it requires full knowledge of the system’s state. However, one may instead use the stochastic phase to reconstruct an approximate state lying on the center manifold. This approximation is expected to remain close to the true state of the system, provided the noise is small and the system stays within the region of stability.

This approximation has two important consequences. First, it allows us to track the system’s evolution using a reduced process whose dimension matches that of the center manifold. Second, in systems with translation invariance—as in our case—this reduction yields an effective Itô process describing the phase dynamics of the wave’s position. Together, these ideas form the basis for the dimensional reduction central to \textbf{isochronal phase reduction} method.

Such a reduction has been rigorously justified by Katzenberger  for systems with finite-dimensional state spaces \cite{Katzenberger}. Parsons and Rogers provide an in-depth tutorial on how to compute such reductions in the finite-dimensional setting \cite{Parsons}. Moreover, there are biologically motivated examples \cite{jones,constable} in which applying this framework to stochastically forced systems has revealed unexpected dynamical behavior. Notably, \cite{jones} studies a spatially extended version of their system using an ad hoc definition of stochastic phase, estimated via Monte Carlo simulation. The method developed in this paper may help formalize such analyses by providing a clear and principled definition of stochastic phase, along with an efficient method for its computation.

For certain classes of infinite-dimensional systems, the validity of this approximation has also been recently established; see the work of Adams and MacLaurin \cite{MacLaurin}. However, to our knowledge, no explicit calculation, formulation, or testing of a method has yet been carried out for specific systems involving pulse propagation. Our contribution is to lay out a concrete procedure for computing the coefficients of the reduced process that approximates the stochastic phase using the isochronal map, in a setting that is broadly applicable and different than previous works. 

We conclude by saying that work on methods for computing a stochastic wave speed has grown quite deep over recent years. The relatively recent trend began with formal and ad hoc methods; see for example work on stochastic neural fields \cite{Bressloff} or the work using so-called collective coordinates \cite{Cartwright} for the Nagumo equation. There is also a long list of work regarding variational stochastic phase by which the phase is constructed via the minimization of some semi-norm measuring the error between the stochastic solution and center manifold. See \cite{Kuehn} for an example of the variational approach for the Fitzhugh-Nagumo equation where the noise is additive and only in the equation for $u$. MacLaurin and Inglis as well as others have used the variational phase approach in neural field models \cite{Inglis, Stannat, LAng}. See \cite{MacLaurin} for a theoretical comparison between variational and isochronal techniques.

Our viewpoint is not to argue for one technique over the others in general, but rather to provide a concrete application of the isochronal phase reduction in an infinite dimensional setting whereas yet none have existed. In fact, the results of Hamster and Hupkes as well as some of the variational phase results are more rigorous and fine-grained in their setting than what we can achieve in our setting with the isochronal phase reduction alone, but the technique presented here is probably more broadly applicable.

In general, one only needs to find a local coordinate system where one coordinate parametrizes the center manifold. In our setting, and what we present in the sequel, we actually find such a global coordinate system, because of the nature of the kinematic equations, but this global-ness is unnecessary. One can, for example, use the eigenvector in the null-space of the adjoint of the linearized operator about the wave solution to co-ordinatize the center manifold as seen in some of the previously mentioned works. 

Part of our reason for choosing the kinematic equation is because it exhibits less regularity--it is piecewise defined, has piecewise smooth solutions, utilizes moving boundaries, does not have an immediately contractive semi-group--than many of the equations used in studying stochastic wave speed. Yet, with a small coordinate transformation, the isochronal phase reduction proves robust against these irregularities, accurate, and relatively simple to implement. Although many of the details in our implementation are specific to the kinematic equations as derived from the FitzHugh–Nagumo model, the computational framework will extend to other PDEs supporting traveling wave solutions and potentially to systems with higher dimensional center manifolds as well. 
 
\section{The Model}
We study a kinematic model that arises as the singular limit $\epsilon \to 0$ of pulse solutions to certain FitzHugh–Nagumo-type equations. The system is defined on a periodic spatial domain of length $L$, which we identify with the quotient space $\mathbb{R}/(L\mathbb{Z})$. Each point $x \in \mathbb{R}/(L\mathbb{Z})$ will be represented by its real-valued coordinate in the interval [0, L).

Let $w(t,x):[0,\infty) \times \mathbb{R}/(L\mathbb{Z}) \to \mathbb{R}$ be continuous and piecewise smooth, and let $x_+(t),x_-(t):[0,\infty) \to \mathbb{R}/(L\mathbb{Z})$ be continuously varying interface positions. Given initial data $(w(0,x),x_+(0),x_-(0))$, we say the triple $(w(t,x),x_+(t),x_-(t))$ satisfies the kinematic equations if 
\begin{equation}
\begin{aligned}
    \partial_t w  &= \begin{cases}
    g_E(w) & x  \in E_{(x_-,x_+)}  \\
    g_R(w) & x \in R_{(x_-,x_+)}
    \end{cases}
    \\
    \dfrac{dx_{\pm}}{dt} &= \pm c(w(t,x_\pm))
\end{aligned}
\label{PDE:Kin_Det}
\end{equation}
where $g_{E},g_{R}:\mathbb{R} \to \mathbb{R}$ and $c:\mathbb{R}\to\mathbb{R}$ are smooth functions.

The sets $E_{(x_-,x_+)}$ and $R_{(x_-,x_+)}$ are defined as the two open, connected intervals in $\mathbb{R}/(L\mathbb{Z})$  whose union covers the domain (up to full Lebesgue measure), and whose common boundary is given by $\{x_-, x_+\}$. That is, $E$ and $R$ are disjoint and together cover the domain (modulo the measure-zero set $\{x_-, x_+\}$), with transitions occurring at $x_-(t)$ and $x_+(t)$. We assume that the choice of $g_E$, $g_R$, and $c$ ensures that $w(t,x)$ remains continuous in $x$ for each fixed $t \geq 0$, so that the evaluation $w(t, x_\pm)$ is always well defined.

 Motivated by applications such as ion channel noise, we study the stochastic version of the kinematic model \eqref{PDE:Kin_Det}, in which noise appears only in the equation for $w$. We consider noise that is white in time, spatially correlated, and modeled in the Itô sense. Additionally, we assume that the noise is \textbf{translationally invariant}, perhaps reflecting the spatial symmetry of an underlying microscopic model. The resulting SPDE system is given by:
\begin{equation}
\begin{aligned}
    dw  &= \sigma h(w)dW + \begin{cases}
    g_E(w)dt & x  \in E_{(x_-,x_+)}  \\
    g_R(w)dt & x \in R_{(x_-,x_+)}
    \end{cases}  
    \\
    dx_{\pm} &= \pm c(w(t,x_\pm))dt
\end{aligned}
\label{SPDE:Kin_Sto}
\end{equation}
where $\sigma > 0$ denotes the noise amplitude, and $W(t,x)$ is a spatially correlated, temporally white noise process on the periodic domain.

\begin{remark}
Throughout the paper, we will occasionally refer to “the deterministic version of the equation,” by which we mean the PDE obtained by formally omitting the noise terms. For instance, the deterministic version of  \eqref{SPDE:Kin_Sto} is the system  \eqref{PDE:Kin_Det}.
\end{remark}
\vspace{2pt}
Since the domain is compact and the noise is translationally invariant, we represent the noise process $W(t,x)$ via a Fourier expansion:

\begin{equation}
    dW(t,x) = a_0dW_0(t) +\sum_{n = 1}^{\infty}a_{n}\sin\left(\frac{2\pi n x}{L}\right)dW_{2n-1}(t)+a_{n}\cos\left(\frac{2\pi n x}{L}\right)dW_{2n}(t)
    \label{E:Noise_expansion}
\end{equation}
where $\{W_n(t)\}_n$ are independent, scalar Brownian motions in time. For such an expansion, taking the coefficient of sine and cosine to be the same for each $n$ is equivalent to taking the noise to be translationally invariant. Letting $\mathbb{E}$ denote expectation, the noise has correlation structure$$d\mathbb{E}(W(t,x),W(t,y)) = C(x-y)dt$$ where $C$ is an even function. Specifically $C$ takes the form
\begin{equation}
    C(x-y) = a_0^2 + \sum_{n=1}^{\infty}a_{n}^2\cos\left(\frac{2\pi n (x-y)}{L}\right),
\label{Correlation}
\end{equation}
so different choices of correlation functions result in different Fourier coefficients $a_n^2$.

We choose the coefficients that satisfy the normalization condition \begin{equation}
    1 = \sum_{n=0}^{\infty}a_n^2,
    \label{Normalization}
\end{equation}
so that the overall noise strength is controlled by the prefactor $\sigma$ in \eqref{SPDE:Kin_Sto}. In what follows, we focus on the regime of small $\sigma$, in which our asymptotic and numerical results are valid.

 Under choices of $g_{E}, g_{R}$ and $c$ considered here, there exists a traveling pulse solution to \eqref{PDE:Kin_Det} 
 of the form
 $$(w_*,x_+,x_-) =(w_*(x-c_0t), x_{+,*}+c_0 t, x_{-,*}+c_0 t))$$
where $w_*$ is a piecewise smooth, continuous function defined on the torus $\mathbb{R}/(L\mathbb{Z})$ and satisfies periodic boundary conditions. This solution represents a pulse profile translating with constant speed $c_0$, with interfaces $x_{+}$ and $x_{-}$ tracking the excited and resting regions. 
To analyze fluctuations around this wave, we rewrite the SPDE \eqref{SPDE:Kin_Sto} in the comoving frame via the coordinate transformation $x \to x+ c_0 t$. Since the noise is translationally invariant, this change does not alter the statistical structure of the noise, allowing us to write a statistically equivalent SPDE:
  \begin{equation}
\begin{aligned}
    dw  &= \sigma h(w)dW+ c_0\partial_x wdt + \begin{cases}
    g_E(w)dt & x  \in E_{(x_-,x_+)}  \\
    g_R(w)dt & x \in R_{(x_-,x_+)}
    \end{cases} 
    \\
    dx_{\pm} &= \left(\pm c(w(t,x_\pm))-c_0\right)dt.
\end{aligned}
\label{PDE:Kin_Sto_co_move}
\end{equation}
To interpret the derivative $\partial_x w$ at the moving interfaces $x_{\pm}(t)$, we take the \textbf{forward derivative:}
\[\partial_x^+w(x_{\pm},t) = \lim_{h \to 0}\dfrac{w(x_++h,t)-w(x_+,t)}{h}.\]
Although the singular limit analysis does not prescribe a one-sided derivative explicitly, the resulting condition $\partial_x w(x_\pm,t) = 0$ (see \cite{Keener1}) implies that the spatial derivative of $w$ must vanish near the interfaces. For numerically computed pulse profiles that have $c_0>0$, this condition is only satisfied when taking a forward derivative at the interface. Without loss of generality, we focus only on the case where $c_0>0$ here. This supports our choice of forward differencing as consistent with the asymptotic structure of the solution.

The deterministic version of the system \eqref{PDE:Kin_Sto_co_move} possesses a one dimensional manifold of fixed points given by
\begin{equation} \mathcal{S} = \{(w_*(\cdot -s),(x_{+,*} +s)_L, (x_{-,*}+s)_L \ | s \in \mathbb{R}\}
\label{Center_m}
\end{equation}
where $(x)_L:= x \mod L.$ 
$\mathcal{S}$ corresponds to translations of the traveling pulse solution. Our goal is to understand how solutions drift stochastically along this manifold. However, stochastic perturbations may cause the interfaces $x_+$ and $x_-$ to fluctuate differently, even if their average motion (i.e., drift) is the same. To separate this collective motion from relative fluctuations, it is useful to introduce a time-dependent coordinate system that keeps the excited and resting regions fixed in space.

\subsection{Change of Coordinates: Fixing the Boundaries}
\label{S:Serious_Change_coordinates}
We define a piecewise, time dependent coordinate transformation that maps the moving interfaces $x_\pm(t)$ to fixed reference points $x_{\pm,*}$ chosen s.t. $(w_*,x_{+,*},x_{-,*})$ is the stable pulse solution of the deterministic version of \eqref{PDE:Kin_Sto_co_move}. The points $x_{\pm,*}$ are well defined up to a translation of their mean value which may be set arbitrarily in the interval $[0,L).$  

Specifically, we introduce the map
\begin{equation}
    \xi(x,x_+,x_-) = 
    \begin{cases} 
    \left(\dfrac{(x-x_-)_L(x_{+,*} - x_{-,*})_L}{(x_+-x_-)_L} +x_{-,*}\right)_L & x \in E_{(x_-,x_+)} \\[.5cm]
    \left(\dfrac{(x-x_+)_L(x_{-,*} - x_{+,*})_L}{(x_--x_+)_L} +x_{+,*} \right)_L & x \in R_{(x_-,x_+)}. 
    \end{cases}
\end{equation}
Here, the notation $(\cdot)_L$ denotes modulo-$L$ arithmetic on the torus $\mathbb{R}/(L\mathbb{Z})$, so that $\xi \in [0,L)$. This map satisfies
$$\  \xi(x_+) = x_{+,*}, \quad  \xi(x_-)= x_{-,*}$$ regardless of which branch is used.  Moreover, the restrictions $$\xi: E_{(x_-,x_+)} \to E_{(x_{-,*},x_{+,*})} \quad \text{and} \quad  \xi: R_{(x_-,x_+)} \to R_{(x_{-,*},x_{+,*})} $$ are continuous bijections as long as $x_+ \neq x_-.$ We are not concerned with the degenerate case $x_+ = x_-$ as such a configuration would correspond to the extinction of the pulse.

We denote the inverse transformation by  $x(\xi,x_+,x_-)$ which is given by
\begin{equation}
    x(\xi,x_+,x_-) = \begin{cases}

    \left(\dfrac{(\xi-x_{-,*})_L(x_+-x_-)_L}{(x_{+,*}-x_{-,*})_L} +x_-\right)_L & \xi \in E_{(x_{-,*},x_{+,*})}\\[.5cm]
    \left(\dfrac{(\xi-x_{+,*})_L(x_--x_+)_L}{(x_{-,*}-x_{+,*})_L} +x_+\right)_L & \xi \in R_{(x_{-,*},x_{+,*})}.   
    \end{cases}
\end{equation}
Since $x_+$ and $x_-$ are stochastic processes above, then the inverse transformation, denoted by $x$ for notational convenience, is also a stochastic process.
However, because $x$ depends linearly on $(x_+,x_-)$ and the dynamics $x_+$ and $x_-$ do not contain explicit stochastic terms, no Itô correction or explicit stochastic terms arise: 
\begin{equation}
    dx = \begin{cases}
    \left(\dfrac{(\xi-x_{-,*})_L(c(w(x_+)+c(w(x_-))}{(x_{+,*}-x_{-,*})_L} -c(w(x_-)) -c_0\right)dt & \xi \in E_{(x_{-,*},x_{+,*})}\\[.5cm]
    \left(\dfrac{(\xi-x_{+,*})_L(-c(w(x_-))-c(w(x_+)))}{(x_{-,*}-x_{+,*})_L} +c(w(x_+)) -c_0\right)dt & \xi \in R_{(x_{-,*},x_{+,*})}.
    \end{cases}
    \label{E:dx}
\end{equation}
Furthermore, any Itô correction vanishes when applying the chain rule to functions of $x$.

Next we compute the spatial derivative of the inverse transformation with respect to $\xi$:
\begin{equation}
    \partial_\xi x(\xi,x_+,x_-) = \begin{cases}

    \dfrac{(x_+-x_-)_L}{(x_{+,*}-x_{-,*})_L}  & \xi \in E_{(x_{-,*},x_{+,*})}\\[.5cm]
    
    \dfrac{(x_--x_+)_L}{(x_{-,*}-x_{+,*})_L} \ & \xi \in R_{(x_{-,*},x_{+,*})}.
    \end{cases}
    \label{E:partial_xi_of_x}
\end{equation}
To carefully perform the coordinate transformation of \eqref{SPDE:Kin_Sto} under the change of variables $x \to \xi$, we define
\[\hat{w}(\xi,t) = w(x(\xi,x_+,x_-),t)\]
 and derive the SPDE satsified by $\hat{w}$. Applying Itô’s formula to the composition of $w$ with the time-dependent spatial map $x(\xi, x_+, x_-)$, we obtain:
 \begin{equation}
 \begin{aligned}
     d\hat{w}(\xi,t) &= d(w(x(\xi,x_+,x_-),t))
     \\
     & =\partial_x w(x(\xi,x_+,x_-), t)dx +dw(x(\xi,x_+,x_-), t).
\end{aligned}
\label{E:Transform_w}
\end{equation}
\begin{remark}
In \eqref{E:Transform_w}, the symbol $d(w)$ in the first line denotes the total stochastic differential of the composite function $w(x(\xi, t), t)$, while $dw$ in the second line refers to the Itô differential acting directly on the second argument $t$ at fixed spatial position.
\end{remark}
\vspace{2pt}
 
The spatial derivative transforms even more directly:
\begin{equation}
    \partial_{\xi} \hat{w}(\xi,t) = \partial_x w(x(\xi,x_+,x_-),t)\partial_\xi x(\xi,x_+,x_-).
    \label{E:Space_Derivative}
\end{equation}

Combining the coordinate-transformed Itô formula \eqref{E:Transform_w}, the spatial derivative \eqref{E:Space_Derivative}, and the original equation \eqref{PDE:Kin_Sto_co_move}, we obtain the SPDE governing $\hat{w}$:
 \begin{equation}
\begin{aligned}
    d\hat{w}  &=\sigma h(\hat{w})d\hat{W}+ \partial_{\xi} \hat{w}\dfrac{c_0dt+dx(\xi,x_+,x_-)}{\partial_\xi x(\xi,x_+,x_-)}  + \begin{cases}
    g_E(\hat{w})dt & \xi  \in E_{(x_{-,*},x_{+,*})}  \\
    g_R(\hat{w})dt & \xi \in R_{(x_{-,*},x_{+,*})}
    \end{cases} 
    \\
    dx_{\pm} &= \left(\pm c(\hat{w}(t,x_{\pm,*}))-c_0\right)dt.
\end{aligned}
\label{PDE:Kin_Sto_co_move_transformed}
\end{equation}
\begin{remark}
    One can readily verify that the triple $(\hat{w}_*, x_{+,*}, x_{-,*})$, corresponding to the traveling pulse in the co-moving frame, is a stationary solution of the determinstic version of \eqref{PDE:Kin_Sto_co_move_transformed}.
\end{remark}
\vspace{2pt}

A first important obersavation is that the noise has been updated in the new coordinates:
\begin{equation}
\begin{aligned}
    d\hat{W}(t,\xi) = dW(t,x(\xi,x_+,x_-)).
    \end{aligned}
\end{equation}

A second important observation is that the transformed SPDE \eqref{PDE:Kin_Sto_co_move_transformed} may be viewed as two SPDEs posed on the disjoint intervals $E_{(x_{-,*},x_{+,*})}$ and $R_{(x_{-,*},x_{+,*})}$ coupled through continuity conditions for $\hat{w}$ at the interfaces $x_{-,*}$and $x_{+,*}$. Although somewhat unwieldy, we now write out the full equation for completeness: 
\begin{equation}
\begin{aligned}
    d\hat{w}  &= \, \sigma h(\hat{w})d\hat{W}
    \\
    \\
    &+\begin{cases} \dfrac{(\xi-x_{-,*})_L(c(\hat{w}(x_{+,*}))+c(\hat{w}(x_{-,*}))}{(x_{+}-x_{-})_L} \partial_{\xi} \hat{w}dt 
      \\[.5cm]
    -\dfrac{(\xi-x_{+,*})_L(c(\hat{w}(x_{+,*}))+c(\hat{w}(x_{-,*}))}{(x_{-}-x_{+})_L}\partial_{\xi} \hat{w}dt 
    \end{cases} 
    \\[.5cm]
    &\begin{array}{ll}
    -\dfrac{c(\hat{w}(x_{-,*}))(x_{+,*}-x_{-,*})_L}{(x_{+}-x_{-})_L}\partial_{\xi} \hat{w}dt 
    \\
    \\
    + \dfrac{c(\hat{w}(x_{+,*}))(x_{-,*}-x_{+,*})_L}{(x_{-}-x_{+})_L}\partial_{\xi} \hat{w}dt
    \end{array}
    \\
    \\
    &\begin{array}{ll}+ \,
    g_E(\hat{w})dt & \quad \xi \in E_{(x_{-,*},x_{+,*})}
    \\
    \\
    +\, g_R(\hat{w})dt & \quad \xi \in R_{(x_{-,*},x_{+,*})}
    \end{array}
    \\[.5cm]
    dx_{\pm} &= \left(\pm c(\hat{w}(t,x_{\pm,*}))-c_0\right)dt.
    \end{aligned}
\end{equation}
The manifold of stationary points for the deterministic version of this equation, originally given in \eqref{Center_m}, now becomes
\[\mathcal{S} = \{(\hat{w}_*, (x_{+,*} +s)_L, (x_{-,*}+s)_L) \ | \ s\in \mathbb{R}\}.\]
In particular, we note that
\begin{equation}
    \pm c(\hat{w}_*(x_{\pm,*}))-c_0 = 0
\end{equation}
since $\hat{w}_*$ is the traveling wave solution with constant speed $c_0$.

Since translation along the manifold is still expressed in terms of the two interfaces $x_+$ and $x_-$, it is convenient to perform one final change of coordinates to isolate the direction along the manifold $\mathcal{S}$.

\subsection{Change of Coordinates: Wave Midpoint and Width}
We introduce one last simpler change of coordinates. We define
\begin{align}
\theta &= \left(x_- +\dfrac{(x_+-x_-)_L}{2}\right)_L
\\
\rho&= \dfrac{(x_+-x_-)_L}{2},
\end{align}
which gives an invertible linear transformation between $(x_+,x_-)$ and $(\theta,\rho).$ Here, $\theta$ is the midpoint and $\rho$ half the width of the excited region $E_{(x_-,x_+)}$. The inverse transformation is
\begin{align}
x_+ &= (\theta +\rho)_L 
\\
x_-&= (\theta -\rho)_L.
\end{align}

All expressions from the previous section may be rewritten in these coordinates:
\begin{equation}
\begin{aligned}
    d\hat{w}  =&\sigma h(\hat{w})d\hat{W}
    \\
    \\&+\begin{cases} \left(  \dfrac{(\xi-x_{-,*})_L(c(\hat{w}(x_{+,*}))+c(\hat{w}(x_{-,*}))}{2\rho}\right)\partial_{\xi} \hat{w}dt 
       \\[.5cm]
     
    \left( -\dfrac{(\xi-x_{+,*})_L(c(\hat{w}(x_{+,*}))+c(\hat{w}(x_{-,*}))}{L-2\rho}\right)\partial_{\xi} \hat{w}dt 
    \end{cases} 
    \\
    \\
    &\begin{array}{ll}
    -\,\dfrac{c(\hat{w}(x_{-,*}))\rho_*}{\rho}\partial_{\xi} \hat{w}dt  +\, g_E(\hat{w})dt & \quad \xi \in E_{(x_{-,*},x_{+,*})}
    \\[.5cm]
    +\,\dfrac{c(\hat{w}(x_{+,*}))(L-2\rho_*)}{L-2\rho}\partial_{\xi} \hat{w}dt +\, g_R(\hat{w})dt & \quad \xi \in R_{(x_{-,*},x_{+,*})}
    \end{array}
    \\
    \\
    d\theta =& \left(\dfrac{c(\hat{w}(x_{+,*}))- c(\hat{w}(x_{-,*}))}{2} -c_0\right)dt
    \\[.5cm]
    d\rho  =& \left(\dfrac{c(\hat{w}(x_{+,*}))+ c(\hat{w}(x_{-,*}))}{2}\right)dt
\end{aligned}
\label{SPDE:Fully_fully_Transformed}
\end{equation}

A key property of the deterministic version of \eqref{SPDE:Fully_fully_Transformed} is that the one-dimensional manifold of fixed points takes the form
\[\mathcal{S} = \{(\hat{w}_*, (\theta_* + s)_L, \rho_*) \ | \ s \in \mathbb{R} \}.\] 
In these coordinates, the variables $\hat{w}$ and $\rho$ evolve independently of $\theta$. Moreover, the time differential for $\theta$ does not require knowledge of its current position reflecting the translation invariance of the original system.

Before computing the stochastic phase, let us examine the transformed noise term more carefully, as it is more subtle than previously acknowledged. A typical term in the expansion \eqref{E:Noise_expansion} takes the form
\[\sin\left(\dfrac{2\pi n x}{L}\right), \] 
which, after applying the transformation $x \to x(\xi,\theta,\rho)$, becomes piecewise due to the definition of $x(\xi, \cdot)$. Explicitly, we have:
\begin{equation}
\begin{aligned}
&\sin\left(\dfrac{2\pi n x(\xi,\theta,\rho)}{L}\right) = \\[.5cm] &\begin{cases}\sin\left(\dfrac{2\pi n (\xi-x_{-,*})_L\rho}{L\rho_*} + \dfrac{2 \pi n(\theta -\rho)}{L}\right)& \xi \in E_{(x_{-,*},x_{+,*})}\\[.5cm]
\sin\left(\dfrac{2\pi n (\xi-x_{+,*})_L(L-2\rho)}{L(L-2\rho_*)} +\dfrac{2 \pi n (\theta +\rho)}{L}\right)& \xi \in R_{
(x_{-,*},x_{+,*})}.\end{cases} 
\end{aligned}
\end{equation}
\begin{remark}
    We have used the $L$-periodicity of $\sin(2 \pi\, \cdot/L)$ to omit the modular arithmetic notation $( \, \cdot \,)_{L}$.
\end{remark}
\vspace{2pt}
Due to the translation invariance of the noise, we may subtract $\theta$ without changing the law, yielding an equivalent process $(\hat{w}, \theta, \rho)$ in distribution. Accordingly, the noise term takes the form
\begin{equation}
\sin\left(\dfrac{2\pi n x(\xi,\rho)}{L}\right) = \begin{cases}\sin\left(\dfrac{2\pi n (\xi-x_{-,*})_L\rho}{L\rho_*} - \dfrac{2 \pi n\rho}{L}\right)& \xi \in E_{(x_{-,*},x_{+,*})}\\[.5cm]
\sin\left(\dfrac{2\pi n (\xi-x_{+,*})_L(L-2\rho)}{L(L-2\rho_*)} +\dfrac{2 \pi n \rho}{L}\right)& \xi \in R_{
(x_{-,*},x_{+,*})}\end{cases}
\end{equation}
where 
\begin{equation}
x(\xi,\rho) = \begin{cases}\dfrac{(\xi-x_{-,*})_L\rho}{\rho_*} - \rho & \xi \in E_{(x_{-,*},x_{+,*})}\\[.5cm]
\dfrac{2\pi n (\xi-x_{+,*})_L(L-2\rho)}{L-2\rho_*} +\rho & \xi \in R_{
(x_{-,*},x_{+,*})}.\end{cases} 
\end{equation}

When computing the reduced Itô equation for $\theta$, we evaluate $x(\xi,\rho)$ at $\rho = \rho_*$, yielding coefficients that depend only on $\hat{w}_*$ and $\rho_*$, and not on $\theta$. This reflects the translation invariance of the system: the update rule for $\theta$ does not require knowledge of its current distribution.

\section{Computing the Stochastic Phase}
We take the isochronal phase as our definition of stochastic phase. Because our coordinate system has been chosen so that the manifold of fixed points is parameterized by the $\theta$ coordinate, defining the isochronal map is relatively straightforward. Considering the deterministic version of \eqref{SPDE:Fully_fully_Transformed}, the isochronal map $\pi$ is defined via the system’s stabilizing dynamics which send initial conditions $(\hat{w}_0, \theta_0, \rho_0)$ \emph{near} the manifold $\mathcal{S}$ to a limiting point $(\hat{w}_*, \theta_\infty, \rho_*)$ \emph{on} $\mathcal{S}$.

Only the second component, $\theta_\infty := \lim_{t \to \infty} \theta(t)$, remains undetermined, and we define this as $\pi_2$. In general, the map $\theta_0 \mapsto \theta_\infty$ may be nonlinear, but near the stable manifold, we only require the first and second derivatives of $\pi_2$ evaluated on $\mathcal{S}$ in order to derive an SDE for $\theta$. Since we are only interested in the second component, we simply write $\pi$ to mean $\pi_2$.

By definition, the isochronal map is given by 
\begin{equation}
    \pi(\hat{w}_0,\theta_0,\rho_0) = \theta_0 + \int_0^{\infty}\left(\dfrac{c(\hat{w}(x_{+,*}))- c(\hat{w}(x_{-,*}))}{2} -c_0\right)dt.
    \label{E:pi}
\end{equation}
To make our notation as transparent as possible, we define 
\begin{equation}
\phi[\hat{w}_0,\rho_0](t,\xi) := \hat{w}(t,\xi)
    \label{Def:phi}
\end{equation} as the trajectory of the $\hat{w}$-component of the solution at time $t$, with initial condition $\hat{w}(0) = \hat{w}_0$ and $\rho(0) = \rho_0$, for the deterministic version of \eqref{SPDE:Fully_fully_Transformed}. We generally suppress the time dependence of $\phi$ below. With this, the map $\pi$ can be rewritten as 
\begin{equation}
    \pi(\hat{w}_0,\theta_0,\rho_0) = \theta_0 + \int_0^{\infty}\left(\dfrac{c(\phi[\hat{w}_0,\rho_0](x_{+,*}))- c(\phi[\hat{w}_0,\rho_0](x_{-,*}))}{2} -c_0\right)dt.
    \label{E:pi2}
\end{equation}

Let us define
\begin{equation}
    \Theta(\hat{w}_0,  \rho_0) := \int_0^{\infty}\left(\dfrac{c(\phi[\hat{w}_0,\rho_0](x_{+,*}))- c(\phi[\hat{w}_0,\rho_0](x_{-,*}))}{2} -c_0\right)dt.
    \label{Def:Theta}
    \end{equation} 

To compute the isochronal stochastic phase we consider the stochastic process $\pi(\hat{w},\theta,\rho)$ and apply the Itô Formula: 

\begin{equation}
    d \pi = \nabla_{\hat{w},\theta,\rho} \pi(\hat{w},\theta,\rho) \{ d\hat{w} , d\theta , d\rho\}+ \dfrac{1}{2} \nabla^2_{\hat{w},\theta,\rho} \pi(\hat{w},\theta,\rho)\{ d\hat{w} , d\theta , d\rho\},
    \label{E:pre_SDE}
\end{equation}
where the notation $\{\cdot\}$ indicates the direction in which each derivative is applied, and derivatives are taken with respect to the $L^2(\mathbb{R}/(L\mathbb{Z}))$ inner product.
\begin{remark}
   The expression  $\nabla^2_{\hat{w},\theta,\rho} \pi(\hat{w},\theta,\rho)\{ d\hat{w} , d\theta , d\rho\} $ represents the evaluation of a quadratic form. Derivatives of $\pi$ with respect to $\hat{w}$ are functional derivatives and must be considered carefully. In particular, this expression should be interpreted as the action of a symmetric bilinear form $B$ evaluated at $ (d\hat{w} , d\theta , d\rho)$ i.e. $B((d\hat{w} , d\theta , d\rho),(d\hat{w} , d\theta , d\rho))$.
\end{remark}
\vspace{2pt}

Now we simplify. Using equations \eqref{E:pi} and \eqref{Def:Theta}, we find 
\begin{equation*}
    d \pi = d\theta+\nabla_{\hat{w},\rho} \Theta(\hat{w},\rho)\{ d\hat{w} ,  d\rho\}+\dfrac{1}{2}\nabla^2_{\hat{w},\rho} \Theta(\hat{w},\rho)\{ d\hat{w} ,  d\rho\}.
\end{equation*}
This expression is not yet closed, as it depends on the full processes $\hat{w}$ and $\rho$. However, since $\hat{w}$ and $\rho$ remain close to $\hat{w}_*$ and $\rho_*$, we approximate by evaluating all terms at this fixed point. This yields:
\begin{equation}
\begin{aligned}
    d \pi_* = \, & d\theta\Big\rvert_{\hat{w}_*,\rho_*}+\nabla_{\hat{w},\rho} \Theta(\hat{w}_*,\rho_*)\{ d\hat{w}\Big\rvert_{\hat{w}_*,\rho_*} ,  d\rho\Big\rvert_{\hat{w}_*,\rho_*}\}
    \\
    &+\dfrac{1}{2}\nabla^2_{\hat{w},\rho} \Theta(\hat{w}_*,\rho_*)\{ d\hat{w}\Big\rvert_{\hat{w}_*,\rho_*} ,  d\rho\Big\rvert_{\hat{w}_*,\rho_*}\}.
    \end{aligned}
    \label{E:pre_SDE2}
\end{equation}
We further simplify using the facts that
\[d\hat{w}\Big\rvert_{\hat{w}_*,\rho_*} = \sigma h(\hat{w}_*)d\hat{W}\Big\rvert_{\rho_*}, \quad d\theta \Big\rvert_{\hat{w}_*\rho_*} = 0, \quad \text{and} \quad d\rho\Big\rvert_{\hat{w}_*\rho_*} = 0.\]
Thus, we arrive at a reduced SDE for $\pi_*$:
\begin{equation}
    d \pi_* = \sigma \nabla_{\hat{w}} \Theta(\hat{w}_*,\rho_*)\{ h(\hat{w}_*)d\hat{W}\}+\dfrac{\sigma^2}{2}\nabla^2_{\hat{w}} \Theta(\hat{w}_*,\rho_*)\{h(\hat{w}_*)d\hat{W}\}.
    \label{E:SDE}
\end{equation}
Now we expand the noise terms in the SDE. The first-order (diffusion) term becomes
\begin{equation}
\begin{aligned}
   &\nabla_{\hat{w}} \Theta(\hat{w}_*,\rho_*)\{ h(\hat{w}_*)d\hat{W}\} = \\ &\sum_{n=0}^\infty a_n\nabla_{\hat{w}} \Theta(\hat{w}_*,\rho_*)\{h(\hat{w}_*)\cos\left(\frac{2\pi n  x(\xi,\rho_*)}{L}\right)\}dW_{2n} \\+&\sum_{n=1}^\infty a_{n}\nabla_{\hat{w}} \Theta(\hat{w}_*,\rho_*)\{h(\hat{w}_*)\sin\left(\frac{2\pi n  x(\xi,\rho_*)}{L}\right)\}dW_{2n-1}.
    \end{aligned}
    \label{E:diffusion_expanded}
\end{equation}
Using independence of the Brownian Motions, the second-order (drift correction) term becomes
\begin{equation}
\begin{aligned}
   &\nabla^2_{\hat{w}} \Theta(\hat{w}_*,\rho_*)\{ h(\hat{w}_*)d\hat{W}\} = \\ &\sum_{n=0}^\infty a_n^2\nabla^2_{\hat{w}} \Theta(\hat{w}_*,\rho_*)\{h(\hat{w}_*)\cos\left(\frac{2\pi n  x(\xi,\rho_*)}{L}\right)\}dt \\+&\sum_{n=1}^\infty a_{n}^2\nabla^2_{\hat{w}} \Theta(\hat{w}_*,\rho_*)\{h(\hat{w}_*)\sin\left(\frac{2\pi n  x(\xi,\rho_*)}{L}\right)\}dt.
    \end{aligned}
    \label{E:drift_expaned}
\end{equation}
We emphasize that the SDE for $\pi_*$ is trivial in the sense that its coefficients depend only on the fixed quantities $\hat{w}_*$, $\rho_*$, and the noise structure. In particular, they do not depend upon $\pi_*$ itself, a result of the translation invariance, and the resulting expression defines an Itô process with time-homogeneous coefficients.

\begin{remark}
    If our original setting was not translationally invariant, then $\theta$ would appear in the coefficients for the SDE for $\pi$. In our approximation, we would then substitute $\pi_*$ for $\theta$, so the SDE for $\pi_*$ would no longer be trivial. The correct coefficients would depend upon the current state of the SDE.
\end{remark}

\subsection{Computing Coefficients of the Stochastic Phase}
From \eqref{E:diffusion_expanded} and \eqref{E:drift_expaned}, we see that computing the diffusion and drift coefficients requires evaluating directional derivatives of the map $\Theta$. We consider a generic direction $v_0:\mathbb{R}/(L\mathbb{Z}) \to \mathbb{R}$. Since $\phi[\hat{w}_*,\rho_*] =\hat{w}_*$, we find that
\begin{equation}
\begin{aligned}
&\nabla_{\hat{w}} \Theta (\hat{w}_*,\rho_*)\{v_0\} = 
\\ &\dfrac{1}{2}\int_0^\infty c'(\hat{w}_*(x_{+,*}))D_{\hat{w}}\phi[\hat{w}_*,\rho_*]\{v_0\}(x_{+,*})dt
\\
 -&\dfrac{1}{2}\int_0^{\infty} c'(\hat{w}_*(x_{-,*}))D_{\hat{w}}\phi[\hat{w}_*,\rho_*]\{v_0\}(x_{-,*})dt. 
\end{aligned}
 \label{Theta_prime}
\end{equation}
The second derivative is given by
\begin{equation}
\begin{aligned}
&\nabla^2_{\hat{w}} \Theta (\hat{w}_*,\rho_*)\{v_0\} = 
\\
 &\dfrac{1}{2}\int_0^\infty c''(\hat{w}_*(x_{+,*}))\left(D_{\hat{w}}\phi[\hat{w}_*,\rho_*]\{v_0\}(x_{+,*})\right)^2dt
\\
 -&\dfrac{1}{2}\int_0^{\infty} c''(\hat{w}_*(x_{-,*}))\left(D_{\hat{w}}\phi[\hat{w}_*,\rho_*]\{v_0\}(x_{-,*})\right)^2dt
\\
+&\dfrac{1}{2} \int_0^{\infty} c'(\hat{w}_*(x_+,*))D^2_{\hat{w}}\phi[\hat{w}_*,\rho_*]\{v_0\}(x_{+,*})dt
\\
-&\dfrac{1}{2} \int_0^{\infty} c'(\hat{w}_*(x_-,*))D^2_{\hat{w}}\phi[\hat{w}_*,\rho_*]\{v_0\}(x_{-,*})dt.
\end{aligned}
\label{theta_primeprime}
\end{equation}
\begin{remark}
    Note that term $D_{\hat{w}}\phi[\hat{w}_*,\rho_*]\{v_0\}$ is a continuous function and can therefore be evaluated point-wise. For any $(t_0,x_0) \in \mathbb{R}^+\times \mathbb{R}/(L\mathbb{Z})$, we write the evaluation as $D_{\hat{w}}\phi[\hat{w}_*,\rho_*]\{v_0\}(t_0,x_0).$ The same holds for $D^2_{\hat{w}}\phi[\hat{w}_*,\rho_*]\{v_0\}$.
\end{remark}
\vspace{2pt}

We thus require first and second derivatives of $\phi$ with respect to $\hat{w}$ in certain directions. For any such direction $v_0$, we compute for the first derivative
\begin{equation} 
v(t):= D_{\hat{w}} \phi [t, \hat{w}_*]\{v_0\}=\lim_{\epsilon \to 0}\dfrac{\phi[t,\hat{w}_*+\epsilon v_0]-\phi[t,\hat{w}_*]}{\epsilon},
\label{E:First_derivative}
\end{equation}
which represents the linear response of the trajectory to a perturbation in the direction $v_0$. The relevant directions are those appearing in the expansion \eqref{E:diffusion_expanded}. While we cannot compute these directional derivatives for infinitely many directions, the decay $a_n \to 0$ allows us to truncate the sum and retain only the dominant contributions. 

Computing \eqref{E:First_derivative} amounts to linearizing \eqref{SPDE:Fully_fully_Transformed} about $\hat{w}_*$. The resulting system is
\begin{equation}
\begin{aligned}
    \dfrac{dv}{dt} =&
    \begin{cases} 

    c_0\partial_{\xi} v -\dfrac{c_0(\rho-\rho_*)}{\rho_*}\partial_{\xi}\hat{w}_*-c_{0,-}'v(x_{-,*})\partial_\xi \hat{w}_*
       \\[.5cm]
    
    c_0\partial_{\xi} v + \dfrac{2c_0(\rho-\rho_*)}{L-2\rho_*}\partial_\xi\hat{w}_*+c_{0,+}'v(x_{+,*})\partial_\xi \hat{w}_*  
    \end{cases} 
    \\
    \\
    &\begin{array}{ll}+\dfrac{(\xi-x_{-,*})_L(c_{0,+}'v(x_{+,*})+c_{0,-}'v(x_{-,*}))}{2\rho_*}\partial_\xi\hat{w}_*   
    \\[.5cm]
    -\dfrac{(\xi-x_{+,*})_L(c_{0,+}'v(x_{+,*})+c_{0,-}'v(x_{-,*}))}{L-2\rho_*}\partial_\xi\hat{w}_*   
    \end{array}
    \\
    \\
    &\begin{array}{ll}
    +\, 
    g_E'(\hat{w}_*)v & \quad \xi \in E_{(x_{-,*},x_{+,*})} 
    \\[.5cm]
    + \,
    g_R'(\hat{w}_*)v & \quad \xi \in R_{(x_{-,*},x_{+,*})}
    \end{array}
    \\
    \\
    \dfrac{d\rho}{dt} & = \left(\dfrac{c_{0,+}'v(x_{+,*})+ c_{0,-}'v(x_{-,*})}{2}\right)
\end{aligned}
\label{PDE:First_Derivative}
\end{equation}
where $c'_{0,-}:= c'(\hat{w}_{*}(x_{-,*}))$  and $c'_{0,+}:= c'(\hat{w}_{*}(x_{+,*}))$. Continuity is enforced at $x_{\pm,*}$ and the initial conditions are $(v(0),\rho(0)) = (v_0,\rho_*).$ 

Before moving on to the second derivative, it is helpful to condense notation. The equation we needed to solve for $(v,\rho)$ is linear in  both $v$ and $\rho-\rho_*$, and so we could rewrite \eqref{PDE:First_Derivative} as 
\begin{equation}
\begin{bmatrix}
    \partial_t v \\
    \frac{d\rho}{dt}
\end{bmatrix} = \begin{bmatrix}
    \mathcal{L}_{1,1} & \mathcal{L}_{1,2} \\
    \mathcal{L}_{2,1} & 0
\end{bmatrix}
\begin{bmatrix}
    v \\ \rho-\rho_*
\end{bmatrix}.
\end{equation}

We now consider the directional second derivative
\begin{equation}\nu(t) := \dfrac{1}{2}D_{\hat{w}}^2\phi[t,\hat{w}_*]\{v_0\}=  \lim_{\epsilon \to 0}\dfrac{\phi[t,\hat{w}_*+\epsilon v_0]-2\phi[t,\hat{w}_*]+\phi[t,\hat{w}_*-\epsilon v_0]}{2\epsilon^2} .
\label{half}\end{equation}

We refer to the Appendix for a more detailed discussion on computing the PDEs for the first and second derivatives. Using \eqref{E:Drift_term}, we obtain the the PDE for $(\nu,\zeta)$:
\begin{equation}
\begin{aligned}
  \partial_t \nu &= \mathcal{L}_{1,1}\nu  +\mathcal{L}_{1,2}(\zeta-\rho_*)
\\
\\
&+\dfrac{1}{2}\begin{cases} \dfrac{(\xi-x_{-,*})_L(c_{0,+}''v(x_{+,*})^2+c_{0,-}''v(x_{-,*})^2)}{2\rho_*}  \partial_{\xi} \hat{w}_* 
       \\[.5cm]
    -\dfrac{(\xi-x_{+,*})_L(c_{0,+}''v(x_{+,*})^2+c_{0,-}''v(x_{-,*})^2)}{L-2\rho_*} \partial_{\xi} \hat{w}_* 
    \end{cases} 
    \\
    \\
    &\begin{array}{ll}+2(\rho_* \partial_\xi v-(\rho - \rho_*) 
 \partial_\xi \hat{w}_*)\dfrac{(\xi -x_{-,*})_L(c_{0,+}'v(x_{+,*})+c_{0,-}'v(x_{-,*}))}{2\rho_*^2}  
    \\[.5cm]
    -2((L-2\rho_*)\partial_\xi v + 2(\rho-\rho_*)\partial_{\xi}\hat{w}_*)\dfrac{(\xi-x_{+,*})_L(c_{0,+}'v(x_{+,*})+c_{0,-}'v(x_{-,*}))}{(L-2\rho_*)^2} 
    \end{array}
    \\
    \\
    &\begin{array}{ll} -2c_{0,-}'v(x_{-,*})\partial_\xi v -c_{0,-}''v(x_-,*)^2  \partial_\xi \hat{w}_*  
    \\[.5cm]
    +2c_{0,+}'v(x_{+,*})\partial_\xi v + c_{0,+}''v(x_+,*)^2 \partial_\xi \hat{w}_* 
    \end{array}
    \\
    \\
    &\begin{array}{ll}
    +2c_{0,-}'\dfrac{(\rho-\rho_*)}{\rho_*} v(x_{-,*})\partial_{\xi}\hat{w}_* - 2c_0\dfrac{(\rho-\rho_*)}{\rho_*}\partial_\xi v + 2c_0 \dfrac{(\rho-\rho_*)^2}{\rho_*^2}\partial_{\xi}\hat{w}_* 
    \\[.5cm]
    + 4c_{0,+}'\dfrac{(\rho-\rho_*)}{L-2\rho_*}v(x_{+,*})\partial_\xi \hat{w}_*
    +4c_0\dfrac{(\rho-\rho_*)}{L-2\rho_*}\partial_{\xi}v +8c_0\dfrac{(\rho-\rho_*)^2}{(L-2\rho_*)^2}\partial_\xi\hat{w}_*
    \end{array}
    \\
    \\
    &\begin{array}{ll}
    +\,  g_E''(\hat{w}_*)v^2  & \quad \xi \in E_{(x_{-,*},x_{+,*})}
    \\[.5cm]
    +\, g''_R(\hat{w}_*) v^2  & \quad \xi \in R_{(x_{-,*},x_{+,*})}
    \end{array}
    \\
    \\
    d\zeta  =& \mathcal{L}_{2,1} \nu + \left(\dfrac{c_{0,+}''v(x_{+,*})^2+ c_{0,-}''v(x_{-,*})^2}{4}\right)dt.
\label{PDE:Second_Derivative}
\end{aligned}
\end{equation}
where $c''_{0,-}:= c''(\hat{w}_{*}(x_{-,*}))$  and $c''_{0,+}:= c''(\hat{w}_{*}(x_{+,*}))$. The initial conditions are $(\nu(0),\zeta(0)) = (0,\rho_*)$ and we enforce continuity at the interfaces $x_{\pm,*}$. 

In summary, the second directional derivative corresponds to $2\nu$ of the pair $(\nu, \zeta)$ which satisfies a linear equation with a relatively large but explicit forcing term constructed from the solution to the first variation $(v, \rho)$. While the expression is lengthy, it consists only of algebraic operations and spatial derivatives and is straightforward to compute in a numerical solver.

\subsection{The Algorithm}
\label{S:Algorithm}
Having derived expressions for the first and second directional derivatives of the flow map $\phi[ \hat{w},\rho]$, we now describe the numerical procedure used to compute the effective drift and diffusion coefficients for the reduced stochastic phase equation. The method proceeds in several steps:
   \begin{itemize}
    \item \textbf{1. Direction selection:} Choose a finite number of directions $v_k$ corresponding to the most significant terms in the expansion of the noise.

    \textbf{For each chosen direction $v_k$, proceed as follows:}

    \item \textbf{2. First variation:} Solve the linearized equation \eqref{PDE:First_Derivative} for the first variation $v(t)$ and $\rho(t)$, using $v_k$ as the initial condition for $v$ and $\rho_*$ as the initial condition for $\rho$.

    \item \textbf{3. Second variation:} Solve the second variation equation \eqref{PDE:Second_Derivative} for $\nu(t)$ and $\zeta(t)$ using $0$ as the initial condition for $\nu$ and $\rho_*$ as the initial condition for $\zeta$. Ensure $v(t)$ and $\rho(t)$ from step 2 are plugged into the source term.

    \item \textbf{4. Integral evaluation:} Compute the integrals in \eqref{Theta_prime} and \eqref{theta_primeprime} using the solutions obtained in steps (2) and (3).

    \textbf{Repeat steps (2)–(4) for all chosen directions.}

    \item \textbf{5. Coefficient assembly:} Assemble the contributions to the drift and diffusion using the expressions in \eqref{E:diffusion_expanded} and \eqref{E:drift_expaned}.
\end{itemize}
We now elaborate on each of these steps in more detail. In particular, step (1) involves a truncation of the noise expansion, which we explain here.

\subsubsection{Step 1: Direction selection}
 Recall the normalization condition
 \begin{equation*}
     1 := \sum_{k=0}^{\infty}a^2_k,
 \end{equation*}
 which ensures $\sigma$ accurately reflects the total energy of the noise. When truncating the expansion, we require that the omitted terms contribute less than the leading-order terms in the SDE coefficients \eqref{SDE:Final}. To this end, we define $n$ to be the smallest natural numbers such that
 \begin{equation*}
\sum_{k>n}a_k \ll \sigma ^2 \quad \text{and} \quad \sum_{k>n}a^2_k \ll\sigma^2. 
 \end{equation*}
 
 After truncation, the directions have been chosen. For $k \in \{0,1,\ldots n_1\}$, the direction vectors are given by 
 \begin{equation}
    \{h(\hat{w}_*)\cos\left(\dfrac{2\pi kx(\xi,\rho_*)}{L}\right)\}_{k=0}^{n}, \quad  \{h(\hat{w}_*)\sin\left(\dfrac{2\pi kx(\xi,\rho_*)}{L}\right)\}_{k=1}^{n}.
 \end{equation}
For each pair of directions, corresponding to index $k$, we perform steps (2)-(4).
 \subsubsection{Step 2: First variation}
We are required to compute 
 \begin{equation}
 \begin{aligned}
\nabla_{\hat{w}}\Theta(\hat{w_*},\rho_*)\{h(\hat{w}_*)\cos\left(\dfrac{2\pi kx(\xi,\rho_*)}{L}\right)\} 
\\
\nabla_{\hat{w}}\Theta(\hat{w_*},\rho_*)\{h(\hat{w}_*)\sin\left(\dfrac{2\pi kx(\xi,\rho_*)}{L}\right)\}.
\end{aligned}
 \end{equation}
To do this, we compute $(v^{(c)}_k,\rho_{k}^{(c)})$ and $(v^{(s)}_k,\rho_{k}^{(s)})$, the solutions to the PDE \eqref{PDE:First_Derivative} with initial conditions 
\begin{align*}
    (v^{(c)}_k(\xi, 0),\rho^{(c)}_k(0))= \left(h(\hat{w}_*)\cos\left(\dfrac{2\pi kx(\xi,\rho_*)}{L}\right),\rho_*\right)
    \\
    (v^{(s)}_k(\xi, 0),\rho^{(s)}_k(0))= \left(h(\hat{w}_*)\sin\left(\dfrac{2\pi kx(\xi,\rho_*)}{L}\right),\rho_*\right).
\end{align*}
We solve each PDE over a truncated time domain $[0,T_s]$ such that sufficient decay of the solution occurs. We formally specify the required decay in step (4), but in practice it is easy to choose a sufficiently large $T_s$ and typically $T_s< 100$.
\subsubsection{Step 3: Second variation}
We are required to computed
\begin{equation}
\begin{aligned}
\nabla^2_{\hat{w}}\Theta(\hat{w_*},\rho_*)\{h(\hat{w}_*)\cos\left(\dfrac{2\pi kx(\xi,\rho_*)}{L}\right)\} 
\\
\nabla^2_{\hat{w}}\Theta(\hat{w_*},\rho_*)\{h(\hat{w}_*)\sin\left(\dfrac{2\pi kx(\xi,\rho_*)}{L}\right)\}.
\end{aligned}
\end{equation}
To do so we solve the PDE \eqref{PDE:Second_Derivative} with corresponding $(v_k^{(c)},\rho_k^{(c)})$ and $(v_k^{(s)},\rho_{k}^{(s)})$ plugged in for the forcing terms dependent on $v$ and $\rho$. As in the previous step, we solve over the truncated time domain $[0,T_s]$. Initial conditions are given by
\begin{align*}
    (\nu^{(c)}_k(\xi, 0),\zeta^{(c)}_k(0))= \left(0,\rho_*\right)
    \\
    (\nu^{(s)}_k(\xi, 0),\zeta^{(s)}_k(0))= \left(0,\rho_*\right).
\end{align*}
Be careful to remember that $\nu$ corresponds to just half the second derivative as seen in \eqref{half}.
\subsubsection{Step 4: Integral evaluation}
Next we compute the integrals seen in \eqref{Theta_prime} and \eqref{theta_primeprime}. These first pair of integrals are used for the approximate diffusion coefficients in the SDE \eqref{E:SDE}:
\begin{equation}
\begin{aligned}
    S_{2k}(T_s) :=&\; \dfrac{a_k}{2}\int_0^{T_s}c'(\hat{w}(x_{+,*}))v^{(c)}_k(x_{+,*})-c'(\hat{w}(x_{-,*}))v^{(c)}_k(x_{-,*})dt
    \\
    S_{k}(T_s):=& \; \dfrac{a_k}{2}\int_0^{T_s}c'(\hat{w}(x_{+,*}))v^{(s)}_k(x_{-,*})-c'(\hat{w}(x_{-,*}))v^{(s)}_k(x_{-,*})dt.
    \label{computed_var}
\end{aligned}
\end{equation}
This second integral is used in computing the approximate drift coefficient:
\begin{equation}
\begin{aligned}
&R_k(T_s) := a_k^2\Bigg( \dfrac{1}{2}\int_0^{T_s}c''(\hat{w}_*(x_{+,*}))\left(\left(v_k^{(c)}(x_{+,*})\right)^2+\left(v_k^{(s)}(x_{+,*})\right)^2\right)dt
\\
 -&\dfrac{1}{2}\int_0^{T_s} c''(\hat{w}_*(x_{-,*}))\left(\left(v_k^{(c)}(x_{-,*})\right)^2+\left(v_k^{(s)}(x_{-,*})\right)^2\right)dt
\\
+&\int_0^{T_s} c'(\hat{w}_*(x_+,*))(\nu_k^{(c)}(x_{+,*})+\nu_k^{(s)}(x_{+,*}))dt
\\-&\int_0^{T_s} c'(\hat{w}_*(x_{-,*}))(\nu_k^{(c)}(x_{-,*})+\nu_k^{(s)}(x_{-,*}))dt.\Bigg)
\label{Computed_Mean}
\end{aligned}
\end{equation}
A sufficiently large $T_s$ is chosen to ensure that remainder terms are higher order in $\sigma$: 
\[ \lim_{T \to \infty} \left|S_k(T_s)-S_k(T)\right|\ll\sigma, \quad  \lim_{T \to \infty} \left|R_k(T_s)-R_k(T)\right|\ll\sigma.\]

\subsubsection{Step 5: Coefficient assembly}
Finally we assemble the results of steps (2)-(4) into coefficients that are numerical approximations of the coefficients for the SDE in \eqref{E:SDE}. Set
\[ S(T_S) := \dfrac{1}{2}\sum_{k=0}^{n_2}R_k(T_S)\]
Then our computed approximate stochastic phase is given by
\begin{equation}
d\pi_*= \sigma^2 Sdt +\sigma \sum_{k=0}^{n}\left( S_{2k}dW_{2k}+S_{2k+1}dW_{2k+1}\right).
\label{SDE:Final}
\end{equation}
Assuming $\pi_*(0) = \theta(0)$, the variance and mean are given by 
\begin{equation}
\mathbb{E}\left[ \pi_*(t)\right] =\sigma^2 St +\theta(0) \quad \text{and} \quad \text{Var}\left[ \pi_*(t)\right] =t\sigma^2\sum_{k=0}^{n}\left(S_{2k}^2+S_{2k+1}^2\right).
\label{Mean and Var}
\end{equation}
For ease of exposition in the sequel, it helps us to define 
\begin{equation}
    \mu := S \quad \text{and} \quad \nu^2 := \sum_{k=0}^{n}\left(S_{2k}^2+S_{2k+1}^2\right).
    \label{Stats}
\end{equation}
\subsection{Final Remarks on the Algorithm}
We first note that the algorithm presented above is meant to be practical, not provably correct. In particular the truncation rule in step (1) may not be sufficient in pathological cases, because the choice of direction, $v_k$, may make the integrals in the expressions for $S_k$ and $R_k$ large enough to cancel out the smallness of the coefficient $a_k$.

Let us make explicit the dependency on $T_s$ and $n$ of $\mu:=\mu(T_S,n)$ and $\nu:=\nu(T_s,n)$. Then strictly speaking we need to have chosen $T_s$ and $n$ s.t.
\[\lim_{m \to \infty}\lim_{T \to \infty}\left|\mu(T_S,n)-\mu(T,m) \right| < \sigma, \quad \lim_{m \to \infty}\lim_{T \to \infty}\left|\nu(T_S,n)-\mu(T,m) \right| < \sigma.\]
That is the difference between the $\mu$ and $\nu$ computed with and without truncation should only contribute a higher order term, $O(\sigma^3)$, to the coefficients of the SDE.

In practice, we find that in many cases we need to think little about the choice of truncation $T_s$ and $n$ and can often take these to be plenty large enough. Even if $n$ and $T_s$ need to be chosen fairly large, computation of the deterministic PDEs \eqref{PDE:First_Derivative} and \eqref{PDE:Second_Derivative} can be made very efficient using higher order numerical methods. In particular we use the standard 4th order Runge-Kutta method, so we can take a larger time step. In contrast, to obtain the mean drift from \eqref{SPDE:Fully_fully_Transformed} requires sampling. Sampled observables then converge to population values quite slowly at the rate of convergence for the law of large numbers. Further, the existence and accuracy of higher order numerical methods for \eqref{SPDE:Fully_fully_Transformed} is not well established. 

\section{Numerical Experiments and Investigations}

Our goal is to compare the mean wave speed and variance computed using the isochronal map with those obtained by directly sampling from the SPDE \eqref{SPDE:Fully_fully_Transformed}. We focus on the case where, in \eqref{PDE:Fitz}, the excitability parameter is fixed at $\beta = 1$, and the cubic nonlinearity is replaced by a piecewise-defined function: 
\begin{equation}
    f_{\text{cub}}(u) = \begin{cases} -u & u <  \alpha  \\
    1- u & u > \alpha
    \end{cases}
\end{equation}
where $\alpha \in (0, 1/2)$ is a tuning parameter that determines the deterministic wave speed $c_0$, the effective pulse width $\rho_*$, and the traveling wave profile $w_*$ solving \eqref{PDE:Kin_Det}.

In the comoving frame, there are distinct points $x_{-,*}, x_{+,*} \in [0, L)$ such that the solution transitions from the excited region $E$ to the relaxation region $R$ at $x_{-,*}$, and from $R$ back to $E$ at $x_{+,*}$, as $x$ decreases. The pulse solution is continuous, and both $E$ and $R$ are open subsets of the spatial domain.

Solving for $w$ in 
\[ f_\text{cub}(u) - w = 0\]
yields two distinct solutions: an excited state $u=g_E(w)$ and relaxation state $u=g_R(w)$,  given by
\begin{equation} g(w) = 
\begin{cases}
g_E(w) = 1 - (1+\gamma) w & x \in E \\
g_R(w) = -(1+\gamma) w & x \in R.
\end{cases}
\end{equation}
In this case, an analytic expression for $c$ is available. We refer the reader to \cite{Keener1} for further details; for convenience, we reproduce the expression here: 
\begin{equation}
c(w) = \dfrac{1 -2\alpha -2w}{\sqrt{(\alpha +w)(1-\alpha -w)}}.
 \end{equation}
 
We obtain the deterministic PDE for the stationary pulse $w_*$:
\begin{equation}
    0 = c_0 \partial_\xi w_* +g(w_*).
\end{equation}
This equation is translation invariant, so to identify a unique solution, we fix the location of one interface by taking $x_{+,*} = 0$. Since the equation $w_*$ piecewise linear, an explicit solution can be written in terms of the remaining parameters $x_{-,*}$ and $c_0$.  Imposing continuity of $w_*$ at both interfaces $x_{+,*}$ and $x_{-,*}$ leads to a system of two nonlinear equations for the two unknowns. We solve this system using Newton’s method, which allows us to compute $x_{-,*}$ and $c_0$ to high precision.

\begin{figure}
    \centering
    \includegraphics[width=.9\linewidth]{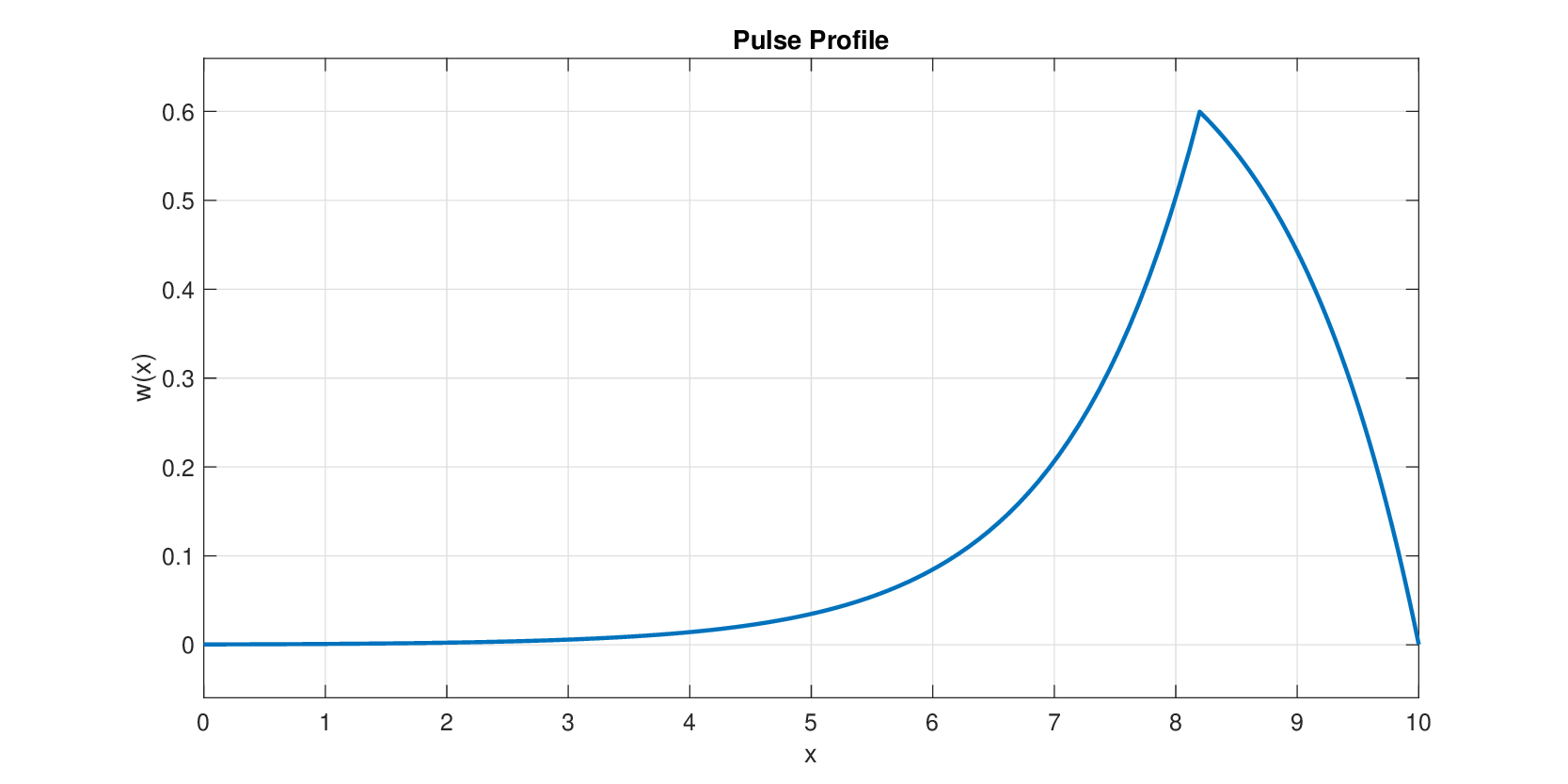}
        \caption{The stationary pulse solution for $L = 10$, $\alpha = 1/5$, and $\gamma = 1/3$. The excited region corresponds to the portion of the pulse where it is decreasing (left to right), while the relaxation region is where it is increasing. The pulse propagates to the right.}
    \label{fig:Stationary}
\end{figure}

A plot of the resulting stationary pulse solution is shown in Figure \ref{fig:Stationary}. As a remark, if the variable $u$ is also of interest, it can be recovered from $w$ via the relations $u = g_E(w)$ in the excited region $E$ and $u = g_R(w)$ in the relaxation region $R$.

We proceed with defining the stochastic setting. As a first pass we consider additive scalar noise by taking  $h \equiv 1$, $a_0 = 1$, and $a_n = 0 $ for $n>0$ corresponding to an additive scalar noise term in the deterministic PDE for $w$. Numerically, we work completely with the transformed equation for $\hat{w}$ \eqref{SPDE:Fully_fully_Transformed}. Note that under this transformation, the scalar additive noise remains unchanged—i.e., it retains the form of a scalar additive noise term.

Our primary goal, at first glance, is to obtain numerical confirmation of the asymptotic theory; specifically, that $\pi_*$, as defined in \eqref{E:SDE}, approximates $\theta$. However, this interpretation requires clarification. The quantity $\pi_*$ is not intended to approximate $\theta$ directly, but rather to approximate the nonlinear projection of the solution onto the center manifold at each time $t$; that is, $\pi(t) = \pi(\hat{w}(t), \theta(t),\rho(t))$ as defined in \eqref{E:pi}. Therefore, while $\theta(t)$ and $\pi(t)$ may not coincide, they are expected to be close, and both are in turn approximated by $\pi_*(t)$. We thus examine the relationship between $\theta$ and $\pi_*$ in our simulations. Of course, we do not use the theoretical expression for $\pi_*$ directly, but rather its computed form given in \eqref{SDE:Final}.

Suppose that $\theta$ indeed evolves according to the reduced stochastic equation \eqref{SDE:Final},  and that we study the drift velocity by computing 
\[ q(T) := \dfrac{\theta(T) - \theta(0)}{T}.\]
Under this assumption, we find that $q(T) \sim \mathcal{N}(\mu\sigma^2,T^{-1}\nu^2\sigma^2)$ where $\mathcal{N}(\mathfrak{m},\mathfrak{d}^2)$ denotes a normal distribution with mean $\mathfrak{m}$ and variance $\mathfrak{d}^2$. Above, $\mu$ and $\nu$ are constants which do not depend on $\sigma$. That is, $q(T)$ is approximately Gaussian with both mean and variance proportional to $\sigma^2$. 

Now suppose we take $N$ realization (i.e. trials) of $q$ and compute their empirical average. Then, 
\[\langle q(T) \rangle \sim \mathcal{N}(\mu\sigma^2, N^{-1}T^{-1}\nu^2 \sigma^2).\] To resolve mean drift $\mu\sigma^2 $ to within tolerance $\epsilon$, we must choose $N$ and $T$ s.t. $N^{-1}T^{-1}\nu^2 \sigma^2 <\epsilon^2$. In particular, to detect a nontrivial drift—i.e., to statistically distinguish it from zero—we must take $\epsilon <\mu \sigma^2.$ This leads to the requirement $N^{-1}T^{-1}<\mu^2 \nu^{-2}\sigma^2$; which can become quite demanding numerically when $\sigma$ is small. Moreover, for certain types of noise—in our experience certain natural multiplicative noise structures—the coefficient $\mu$ itself may be very small, further increasing the computational cost. In such cases, accurately resolving the drift may require taking both $N$ and $T$ to be quite large.

Finally we note here some of the technical details of our simulations. For the full stochastic PDE involving $(\hat{w},\rho,\theta)$ we use the finite difference method with up-winding and the Euler-Maruyama method. We found that the discretization of the deterministic version of \eqref{SPDE:Fully_fully_Transformed} has a slightly different fixed point for $\hat{w}$ and $\rho$ due to the discretization error. The numerical fixed points converge to the true\footnote{By true, we mean the highly accurate solution obtained via Newton’s method.} stable pulse $\hat{w}_*$ and pulse radius $\rho_*$ as the discretization is refined. However, even small mismatches induce a spurious drift in $\theta$, which vanishes in the continuum limit. 

Specifically, we observe that $c(\hat{w}_{\text{numerical}*}(x_{+,*}))\neq -c(\hat{w}_{\text{numerical}*}(x_{{-,*}}))\neq c_0$  and thus $\theta$ drifts with a small constant velocity in the absence of noise. To correct this, we subtract this baseline drift from our measurements. The actual quantity we report is
\[\langle q_{\sigma}(T)\rangle - q_{0}(T).\] 

For the PDEs \eqref{PDE:First_Derivative} and \eqref{PDE:Second_Derivative}, we use finite differences with fourth order Runge-Kutta. Integrals in \eqref{theta_primeprime} are evaluated using the trapezoidal rule. All simulations were performed in MATLAB.

\subsection{Additive Scalar Noise}

Now we finally turn to the results when the noise is an additive scalar. We can approximate $\mu$ by computing $S$ in \eqref{Computed_Mean} without the coefficient $a_0^2$ which is equal to $\sigma^2$ in this example. Similarly we can approximate $\nu$ by computing $S_0$ in \eqref{computed_var} but again without the coefficient $a_0$.  We find $\mu^2/\nu^2$ is roughly $1$. Therefore, we need $NT>>\sigma^{-2}.$ Ideally we shoot for $NT=\sigma^{-3}$. In order to have manageable computational costs, we take $N = T =256$. We also cannot take $\sigma$ too large as the noise can destroy stability. Therefore we take for our experiment 
\begin{equation}
\sigma \in \{\sqrt{2}/16,1/16,\sqrt{2}/32,1/32,\sqrt{2}/64\}.
\label{The Sigma}
\end{equation}

\begin{figure}
    \centering
    \title{\textbf{Results For Scalar Noise}}
\includegraphics[width=0.9\linewidth]{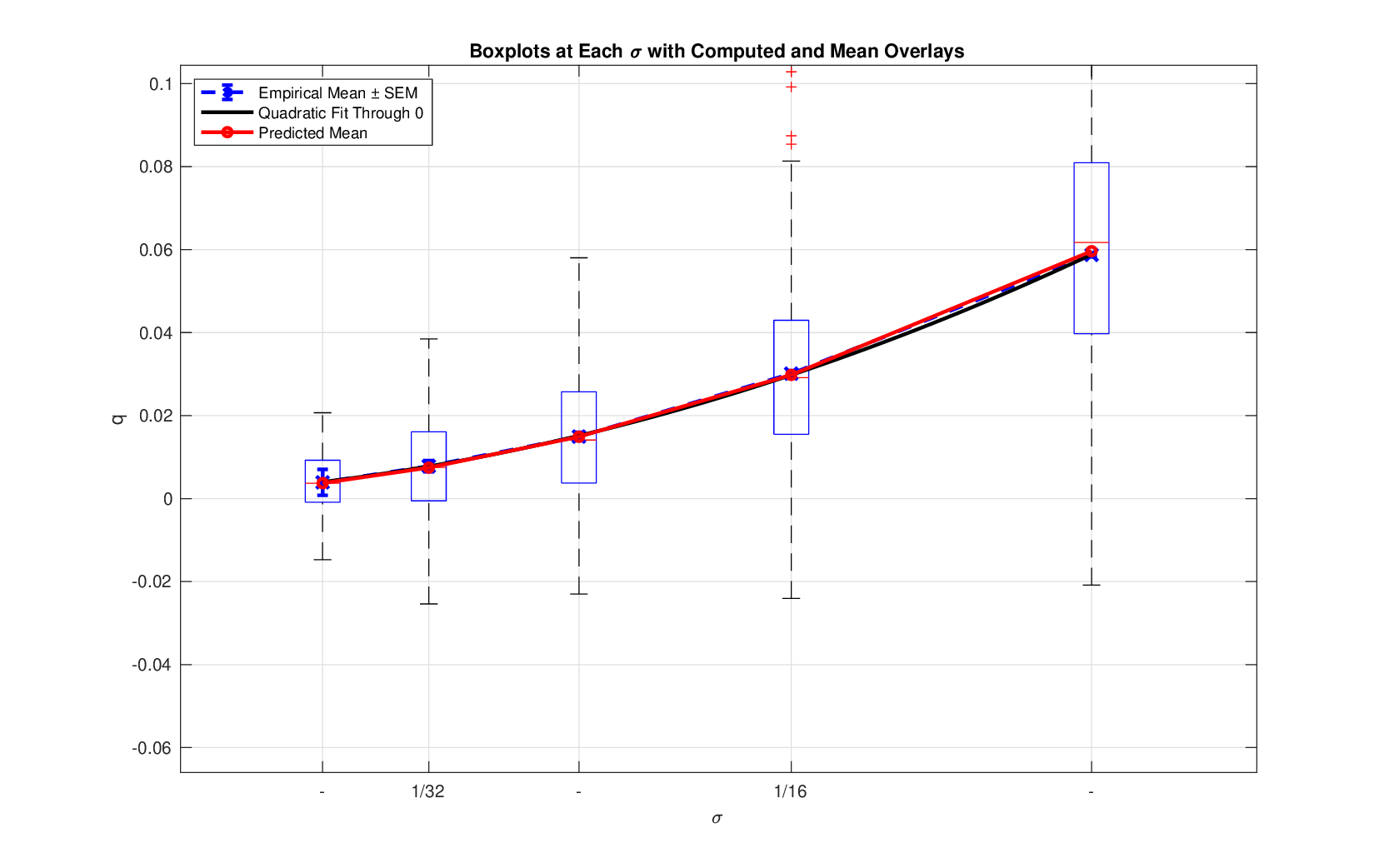}
    \caption{Box plots of the empirical drift are overlayed by empirical means and predicted means. The predicted mean always falls within the standard error of mean (SEM).}
    \label{fig:scalar}
\end{figure}

We refer to $\langle q_\sigma(T)\rangle - q_0(T)$ as the ``Empirical Means'' and we refer to the $S$ in \eqref{Computed_Mean}  as the ``Predicted Means''. We report these in Table \ref{T:Scalar}.
\begin{table}[ht]
\centering
\caption{\textbf{Empirical vs Predicted Means}}
\label{T:Scalar}
\begin{tabular}{|c|c|c|}
\hline
$\sigma$ & \textbf{Empirical Mean} & \textbf{Predicted Mean} \\
\hline
$\tfrac{\sqrt{2}}{64}$ & 0.0037 & 0.0037 \\
$\tfrac{1}{32}$         & 0.0074 & 0.0075 \\
$\tfrac{\sqrt{2}}{32}$ & 0.0148 & 0.0149 \\
$\tfrac{1}{16}$         & 0.0296 & 0.0298 \\
$\tfrac{\sqrt{2}}{16}$ & 0.0591 & 0.0596 \\
\hline
\end{tabular}
\end{table}
We supply box-plots in Figure \ref{fig:scalar} which demonstrate the spread of our observations along with our predicted and empirical means.

\subsection{Additive Noise with Gaussian Correlation}
\begin{figure}[ht]
    \centering
\includegraphics[width=0.9\linewidth]{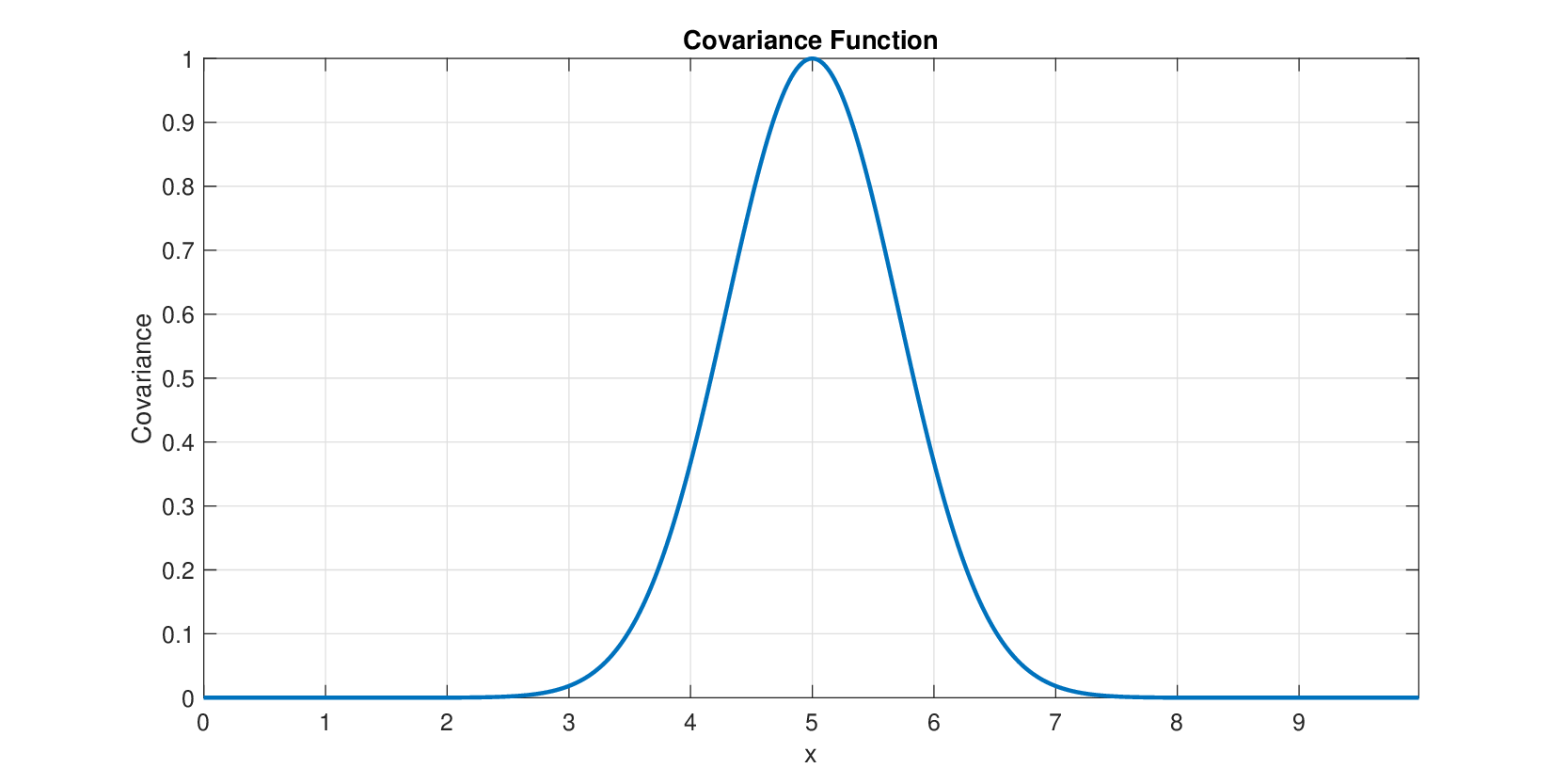}
    \caption{Gaussian covariance function in the variable $x.$ }
    \label{fig:Gauss}
\end{figure}

We now consider a less spatially correlated form of additive noise, in which the coefficients $a_k$ are chosen so that the resulting noise has a correlation function given by the periodization of a Gaussian. The corresponding correlation function is plotted in Figure~\ref{fig:Gauss}.

For simulations of \eqref{SPDE:Fully_fully_Transformed} and for computing the empirical means, we use the first $20$ coefficients of the series in \eqref{E:Noise_expansion} i.e. $\{a_k\}_{k=0}^{19}$. However, to compute the predicted means using the reduced process, we only retain the first $11$ coefficients. The approximate squared values $\{a_k^2\}_{k=0}^{10}$  used in our computation are:
\begin{equation}
\begin{aligned}
   \{a_k^2\}_{k=0}^{10} \approx \{0.17724, \, 0.32118, \, 0.23886, \, 0.14582, \, 0.07308, \, 0.03006, \, \\ 0.01015, \, 0.00281 \, 0.00064 \, 0.00012, \, 0.00002 \}.
    \end{aligned}
\end{equation}

We use slightly different set for $\sigma$ than in \eqref{The Sigma}, as larger values occasionally led to collapse of the wave followed by numerical blow-up of $\theta$, which obscures the data. The set of values we use is:
\[\sigma \in \{1/16,\sqrt{2}/32,1/32,\sqrt{2}/64, 1/64\}. \]

\begin{figure}
    \centering
\title{\textbf{Results for Gaussian Correlation}}    \includegraphics[width=0.9\linewidth]{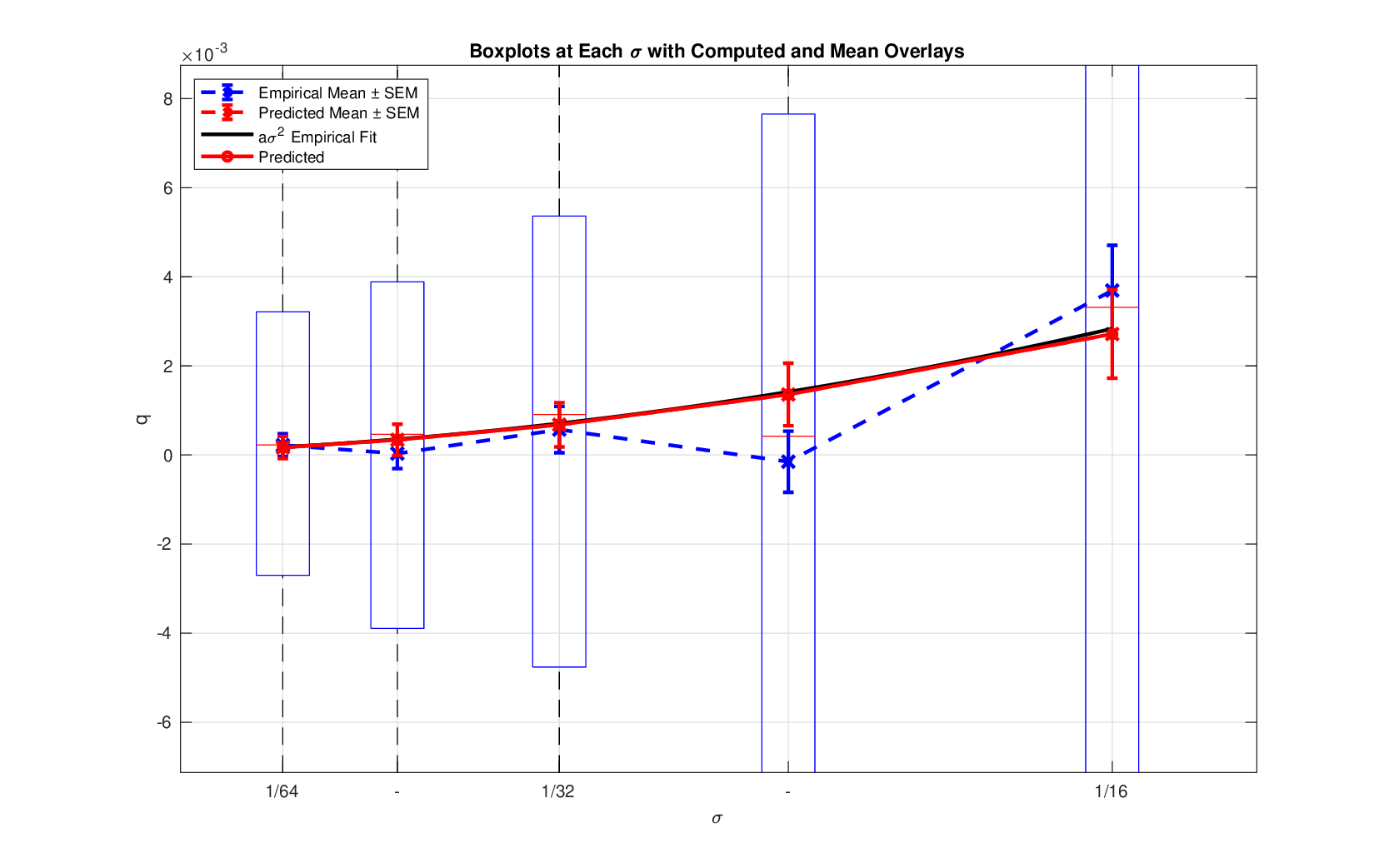}
    \caption{Box plots of the empirical drift are overlayed by empirical means and predicted means. A quadratic is fitted to the empirical means and is in close agreement with the predicted quadratic. SEMs are plotted to indicate statistical imprecision.}
    \label{fig:box}
\end{figure}

We perform an experiment with $T=N=256.$The predicted means, along with the empirical means and medians, are reported in Table \ref{Table:Means}.
\begin{table}[ht]
\centering
\caption{\textbf{Empirical vs Predicted Means}}
\label{Table:Means}
\begin{tabular}{|c|c|c|c|}
\hline
$\sigma$ & \textbf{Empircal Median} &\textbf{Empirical Mean} & \textbf{Predicted Mean} \\
\hline
$\tfrac{1}{64}$ & 0.0002 & 0.0002 & 0.0002 \\
$\tfrac{\sqrt{2}}{64}$ & 0.0005 & 0.0000  & 0.0003  \\
$\tfrac{1}{32}$ & 0.0009 & 0.0006 & 0.0007  \\
$\tfrac{\sqrt{2}}{32}$  & 0.0004 & -0.0002 & 0.0014 \\
$\tfrac{1}{16}$ & 0.0033 & 0.0037 & 0.0027 \\
\hline
\end{tabular}
\end{table}
While these values do not align as closely as in the scalar noise case, the discrepancies are modest. A primary source of inaccuracy arises from the noisiness of the samples relative to the size of the mean. As shown in Figure \ref{fig:box}, which visualizes the results of our experiment, the predicted mean falls within the empirical standard error of the mean (SEM) for all values of $\sigma$ except $\sigma = \sqrt{2}/32$. The SEM, defined as the empirical standard deviation divided by $\sqrt{N}$, provides a measure of the best-case precision attainable from the sample.

To the extent that the empirical means and standard error of the means (SEMs) do not reflect their true population statistics, the discrepancy arises from both numerical and statistical sources of error. The main source of numerical error stems from the time and spatial discretization used to solve \eqref{SPDE:Fully_fully_Transformed}. Specifically, we employed a time-step of $10^{-3}$ and a spatial increment of approximately $0.0615$, corresponding to $605$ grid points.

The primary source of statistical error comes from sample variability and the uncertainty in the sample means. These can be reduced by increasing the number of realizations $N$ and the integration time $T$. However, in our simulations, the time-step, spatial increment, and sample variability are all comparable to or larger than the empirical means themselves. As a result, we should not expect high precision in the first nonzero digit of the empirical means. We attribute the discrepancies between predicted and empirical means primarily to these sources of error.

To support this claim, we fit a quadratic of the form $a\sigma^2$ to the empirical means, choosing $a$ via least squares to minimize the squared error. This form reflects the theoretical expectation that the stochastic drift vanishes when $\sigma = 0$, while assuming a nonzero drift for any nonzero noise.

The fitted curve averages across all sample values of $\sigma$, helping to smooth out numerical and statistical artifacts. As shown in Figure \ref{fig:box}, the fitted quadratic matches the predicted curve nearly exactly, suggesting that the leading-order $\sigma^2$ dependence of the drift is well captured. However, this does not rule out the possibility of higher-order contributions becoming relevant for larger $\sigma$. 

As further evidence supporting the accuracy of our method, we also compare the predicted and empirical standard deviations. Specifically, we compute the predicted standard deviation from the reduced equation as (c.f. \eqref{Mean and Var})
\[\sigma \sqrt{T^{-1}\left(\sum_{k=0}^{11}S_{2k}^2 + S^{2}_{2k+1}\right)},\] and compare this with the empirical standard deviation obtained from sample values of $q$. These quantities are typically larger in magnitude than the empirical means and thus suffer less from relative numerical and statistical error, providing a better basis for comparison. The results are shown in Table \ref{T:dev}.
\begin{table}[ht]
\centering
\caption{\textbf{Empirical vs Predicted Deviations}}
\label{T:dev}
\begin{tabular}{|c|c|c|}
\hline
$\sigma$ & \textbf{Empirical Deviation} & \textbf{Predicted Deviation} \\
\hline
$\tfrac{1}{64}$ &  0.0041 & 0.0040 \\
$\tfrac{\sqrt{2}}{64}$ & 0.0054   & 0.0056 \\
$\tfrac{1}{32}$ & 0.0083  & 0.0079 \\
$\tfrac{\sqrt{2}}{32}$         & 0.0110  & 0.0112 \\
$\tfrac{1}{16}$ & 0.0162  & 0.0159 \\
\hline
\end{tabular}
\end{table}

\subsection{Additive Noise for Single Wave Number}
\begin{figure}[ht]
    \centering
    \title{\textbf{Results for Single Wave Number}}
    \includegraphics[width=0.9\linewidth]{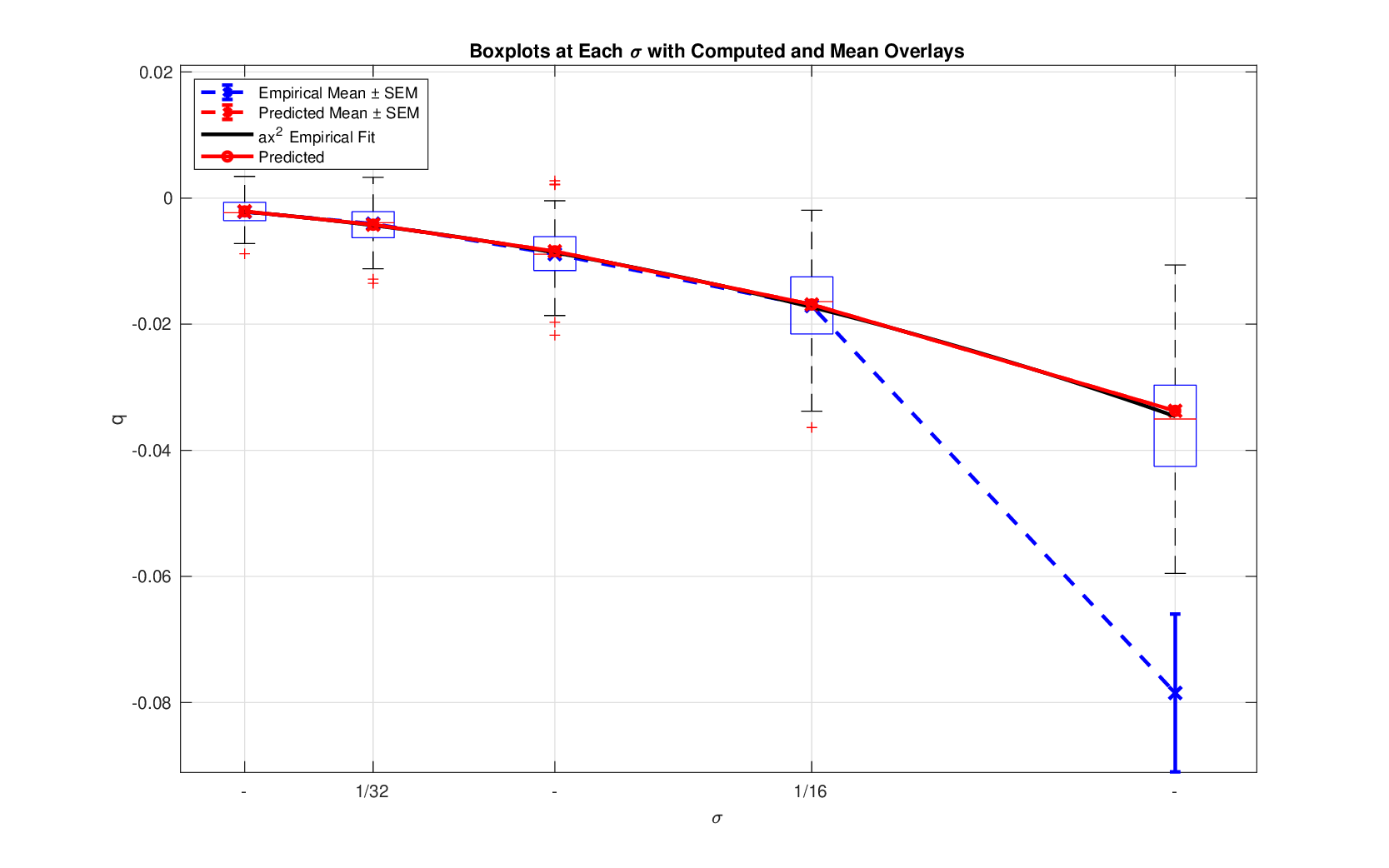}
    \caption{Box plots of the empirical drift are overlayed by empirical means and predicted means. A quadratic is fitted to the empirical \emph{medians} and is in close agreement with the predicted quadratic. SEMs are plotted to indicate statistical imprecision.}
    \label{fig:high_freq}
\end{figure}

We now consider an example where the noise is concentrated at a single wave number:
\begin{equation}
    a_5 = 1 \quad \text{and} \quad a_k =0 \quad \text{otherwise}. 
\end{equation}
This case illustrates how the structure of the noise can qualitatively affect the direction of the stochastic drift. In contrast to the scalar noise example, where the noise increases the average wave speed, the present example shows that the noise instead slows the wave down. This suggests a trend: noise dominated by low-frequency components tends to accelerate the wave, while noise dominated by moderate to high frequencies tends to slow it down. The results are found in Tables \ref{T:Mean_high_freq} and \ref{T:dev_high_freq} and Figure \ref{fig:high_freq}.
\begin{table}[ht]
\centering
\caption{\textbf{Empirical vs Predicted Means}}
\label{T:Mean_high_freq}
\begin{tabular}{|c|c|c|c|}
\hline
 $\sigma$ & \textbf{Empirical Median} & \textbf{Empirical Mean} & \textbf{Predicted Mean} \\
\hline
$\tfrac{\sqrt{2}}{64}$ & -0.00231 & -0.00218 & -0.00206 \\
$\tfrac{1}{32}$         & -0.00390 & -0.00407 & -0.00413 \\
$\tfrac{\sqrt{2}}{32}$ & -0.00890 & -0.00888 & -0.00825 \\
$\tfrac{1}{16}$         & -0.01642 & -0.01711 & -0.01651 \\
$\tfrac{\sqrt{2}}{16}$ & -0.03504 & -0.07859 & -0.03301 \\
\hline
\end{tabular}
\end{table}

\begin{table}[ht]
\centering
\caption{\textbf{ Empirical vs Predicted Deviations}}
\label{T:dev_high_freq}
\begin{tabular}{|c|c|c|}
\hline
 $\sigma$ & \textbf{Empirical Deviation} & \textbf{Predicted Deviation} \\
\hline
$\tfrac{\sqrt{2}}{64}$ &  0.00206 & 0.00211 \\
$\tfrac{1}{32}$         & 0.00313   & 0.00299 \\
$\tfrac{\sqrt{2}}{32}$ & 0.00399  & 0.00423 \\
$\tfrac{1}{16}$         & 0.00615  & 0.00598 \\
$\tfrac{\sqrt{2}}{16}$ & 0.2003  & 0.00845 \\
\hline
\end{tabular}
\end{table}

\subsection{Scalar Multiplicative Noise}

We now consider multiplicative noise where $h \neq 1$. In many biophysical models, it is natural for the intensity of the noise affecting $w$ to depend on the underlying voltage-like variable $u$. For instance, in spatially extended models of stochastic ion channels, the rates of channel opening and closing often vary exponentially with membrane potential \cite{ermentrout2010}.

In our kinematic setting, the voltage variable $u$ is recovered from the relation $f_{\text{cub}}(u) = w$, using the two stable branches of the nonlinearity. In the case where $f_{\text{cub}}$ is defined as a linear function plus a Heaviside step, we recover $u$ via:
\[u(\xi) = \begin{cases}
1- w(\xi) & \xi \in E_{(x_{-,*},x_{+,*})}
\\
- w(\xi) & \xi \in R_{(x_{-,*},x_{+,*})}
    \end{cases}.
\]
We then define the multiplicative noise profile as
\[h(w(\xi)) = e^{u(\xi)}= \begin{cases}
e^{1- w(\xi)} & \xi \in  E_{(x_{-,*},x_{+,*})}
\\
e^{- w(\xi)} & \xi \in R_{(x_{-,*},x_{+,*})}
    \end{cases}.
\]

\begin{figure}[ht]
    \centering
    \title{Results for Multiplicative Noise}
    \includegraphics[width=0.9\linewidth]{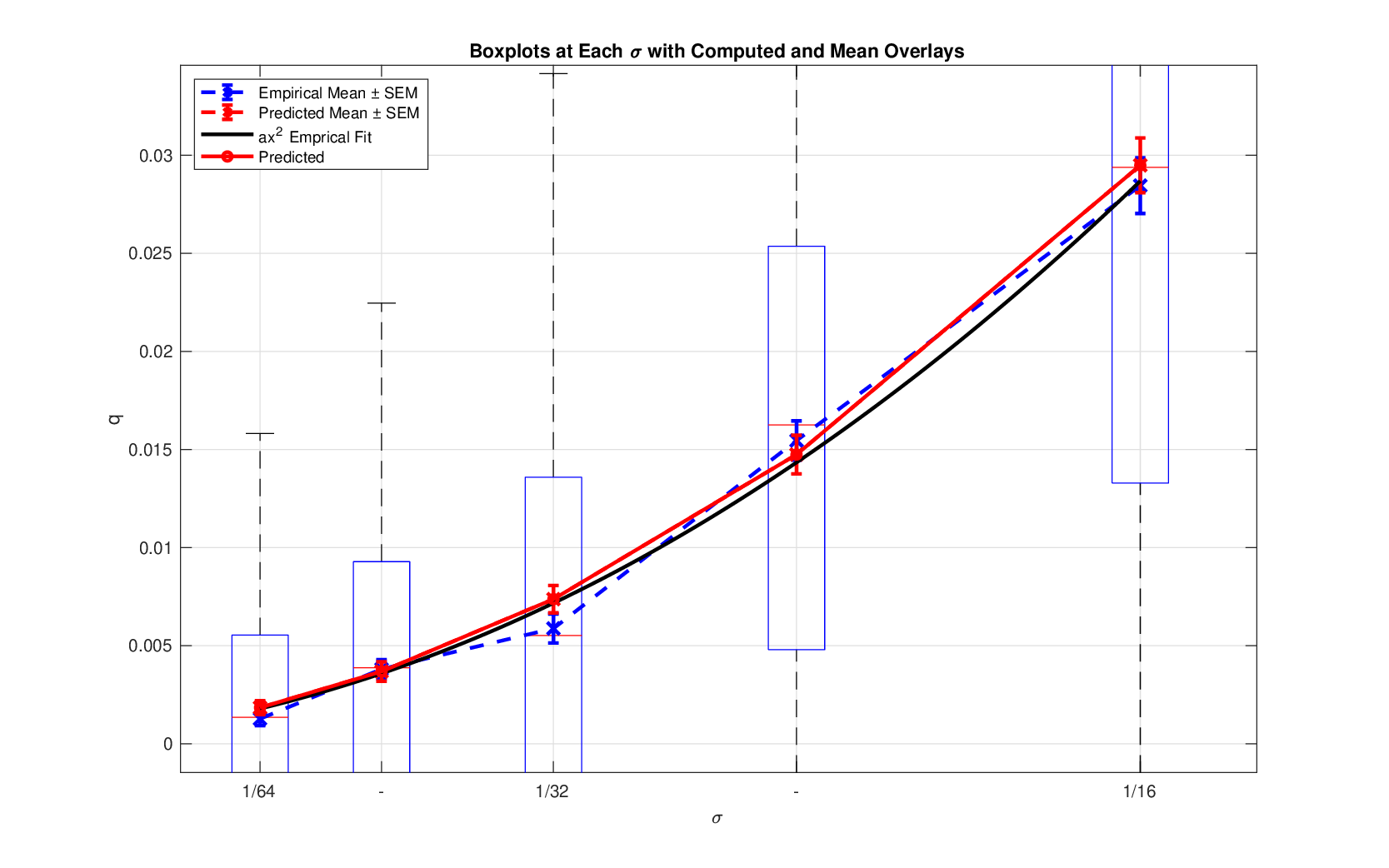}
    \caption{Box plots of the empirical drift are overlayed by empirical means and predicted means. A quadratic is fitted to the empirical means and is in close agreement with the computed quadratic. SEMs are plotted to indicate statistical imprecision.}
    \label{fig:Scalar_Mult}
\end{figure}

This setup encodes the idea that noise intensity grows with increasing depolarization, consistent with biological intuition. The results for this scalar multiplicative noise case are summarized in Tables \ref{Table:Mult_scalar_mean} and \ref{Table:Mult_scalar_dev} and Figure \ref{fig:Scalar_Mult}.
\begin{table}[ht]
\centering
\caption{\textbf{Empirical vs Predicted Means}}
\label{Table:Mult_scalar_mean}
\begin{tabular}{|c|c|c|c|}
\hline
 $\sigma$ & \textbf{Empirical Median}
&\textbf{Empirical Mean}
&\textbf{Predicted Mean} \\
\hline
$\tfrac{1}{64}$ & 0.0013&  0.0013 &0.0018 \\
$\tfrac{\sqrt{2}}{64}$       &  0.0055 & 0.0038   & 0.0037 \\
$\tfrac{1}{32}$ & 0.0055& 0.0059  & 0.0074 \\
$\tfrac{\sqrt{2}}{32}$  &  0.0162      & 0.0155  & 0.0147 \\
$\tfrac{1}{16}$& 0.0294 & 0.0284 & 0.0295 \\
\hline
\end{tabular}
\end{table}

\begin{table}[ht]
\centering
\caption{\textbf{Empirical vs Predicted Deviations}}
\label{Table:Mult_scalar_dev}
\begin{tabular}{|c|c|c|}
\hline
 $\sigma$ & \textbf{Empirical Deviations} &\textbf{Predicted Deviations} \\
\hline
$\tfrac{1}{64}$ &  0.00557 &0.00553 \\
$\tfrac{\sqrt{2}}{64}$         & 0.00768   & 0.00782 \\
$\tfrac{1}{32}$ & 0.01187  & 0.01106 \\
$\tfrac{\sqrt{2}}{32}$         & 0.01600  & 0.01564 \\
$\tfrac{1}{16}$ & 0.02265 & 0.02212 \\
\hline
\end{tabular}
\end{table}
\subsection{Final Remarks on Numerical Experiments}
The results of the previous sections support the idea that the isochronal reduction is very accurate when it comes to predicting the mean drift and drift deviation. Moreover these quantities can be computed efficiently. However, we did not have success with quantifying the error. One expectation might be that the predicted mean and empirical mean agree with each other asymptotically. One hypothesis would be that
\begin{equation}
    \left| \text{Empirical Mean} -\text{Predicted Mean}\right| = O(\sigma^3),
\end{equation}
where the power $3$ is expected since the means themselves are $O(\sigma^2)$. However, we could not find any evidence for this. Indeed our means all agree roughly to the same extent or not at all. See for example Figure \ref{fig:high_freq}. This lack of asymptotic convergence could be due to an insufficient amount of data, which is a result of high computational costs, or it could be due to theoretical considerations.

\subsection{Numerical Investigations: Parameter Exploration}
Finally, we use the algorithm we developed to efficiently explore parameter dependencies of the stochastic drift of the wave. We calculate values independently of $\sigma$ i.e. we compute $\mu$ c.f. \eqref{Stats}. This way the results of this section can be interpreted as valid results for sufficiently small $\sigma$. Although our method is far, far more efficient than Monte-Carlo, we none the less have too many choices to explore everything here. Most of these choices come from the form of $h$ and the coefficients $\{a_k\}_{k=0}^{\infty},$ not to mention we have already picked $f_\text{cub}$. To that end, we fix $h\equiv 1$ and pick the correlation function to be a Gaussian with correlation length $1/\ell^2.$ That is we pick the correlation function $C= C_{\ell,L}$, c.f. \eqref{Correlation}, to be 
\[C_{\ell,L}(x) = \dfrac{\sum_{n\in\mathbb{Z}}e^{-\ell(x+nL)^2}}{N.C.}\]
where $N.C.$ is the normalization constant chosen so \eqref{Normalization} is satisfied.

\begin{figure}[ht]
    \centering
\includegraphics[width=\linewidth]{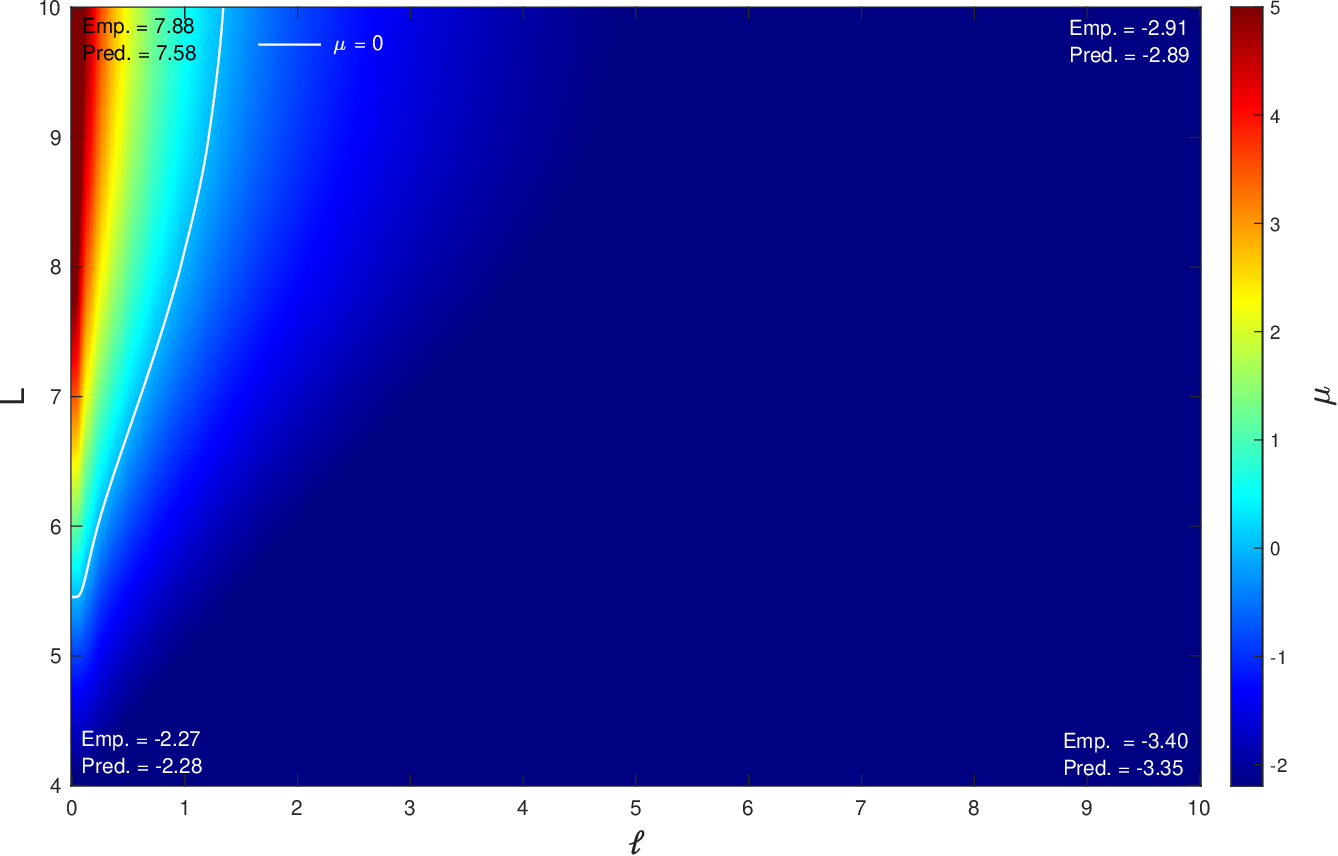}
    \caption{The heatmap shows how relative values of $\mu$ for different values of $\ell$ and $L$. We compare empirical (Emp.) values with predicted (Pred.) values of $\mu$ in each of the four corners of the plot.}
    \label{fig:heat_map}
\end{figure}

We use the method developed in Section \ref{S:Algorithm} to compute various $\mu$ for different values of $\ell$ and $L$. Our findings can be observed in Figure \ref{fig:heat_map}. We find that there exists a small narrow region near $\ell = 0 $ and for large enough domain length for which highly correlated noise causes positive drift of the wave. For less correlated noise and as the correlation length approaches $0$ i.e. the white noise limit, we see negative drift. The limiting value is apparently agnostic to the domain length.

Finally, we test our predictions against empirical observations at values of $\ell$ and $L$ at the four corners of the plot. The corner for which $L=10, \ell = 0$ has already shown in in Table \ref{T:Scalar}. We find the predicted and empirical drift for various values of $\sigma$ for the case of $L=4, \ell = 10$ as can be seen in Table \ref{T:Last}.

\begin{table}[ht]
\centering
\caption{\textbf{Empirical vs Predicted Deviations}}
\label{T:Last}
\begin{tabular}{|c|c|c|}
\hline
 $\sigma$ & \textbf{Empirical Mean} &\textbf{Predicted Mean} \\
\hline
$\tfrac{\sqrt{2}}{64}$         & -0.00134   & -0.00141 \\
$\tfrac{1}{32}$ & -0.00259  & -0.00282 \\
$\tfrac{\sqrt{2}}{32}$         & -0.00520  & -0.00565 \\
$\tfrac{1}{16}$ & -0.011366 & -0.011297 \\
$\tfrac{\sqrt{2}}{16}$ &  -0.0624 &-0.0226 \\
\hline
\end{tabular}
\end{table}



\section*{Declarations}

\textit{Conflicts of Interests:} All four authors declare there is no conflict of interest in this study. 

\noindent\textit{Data Availability:} The numerical data used in this study was generated synthetically using the methods described in the paper. Due to its size and complexity, the raw data has not been made publicly available. However, all necessary details for reproducing the data, including the numerical methods and parameters, are provided in the manuscript. Researchers wishing to replicate the data may follow the outlined methods or contact the authors for further guidance. The corresponding author will share MATLAB code upon request.

\section*{Acknowledgements}
Y.M. was supported by the NSF Materials Research Science and Engineering Center (DMR-2309034). Y.M. and J.M. were supported by the Math+X Award (Proposal Number 234606) from the Simons Foundation.
The authors used OpenAI’s ChatGPT to assist with improving the clarity and style of the text and figures, as well as to facilitate code development. All scientific content, analysis, and interpretations are the authors' own.
\bibliographystyle{siamplain}
\bibliography{sources}

\section*{Appendix}
\label{S:Toy_example}
Here we briefly run through the process of calculating the derivatives of the Isochronal map for a toy example. We consider the system of stochastic differential equations 
\begin{equation}
\begin{aligned}
    du &= f(u,z)dt+ \sigma dW
\\
    dy &= g(u(x_0))dt
\\
    dz & =h(u(x_0))dt
\end{aligned}
\label{E:Toy_SDE}
\end{equation}
where $(u,y,z) \in \mathcal{C}([0,L])\times \mathbb{R} \times \mathbb{R} $. Suppose there is a manifold of attractive fixed points
\[\mathcal{S} = \{(u_*,y,z_*) \ | \ y \in \mathbb{R}\}\] for the deterministic version.

The system captures analogously the relationship between $\hat{w}$, $\theta$, and $\rho$ in \eqref{SPDE:Fully_fully_Transformed}. We may define the isochronal map as before $\pi(u_0,y_0, z_0) =(u_*,y_\infty,z_*)$  where $y_\infty = \lim_{t \to \infty} y(t)$
where 
\[y_\infty = y_0 + \int_0^\infty g(u(x_0))dt.\]
At the end of the day, we want to obtain a reduced process for $y$. Thus we only care about the second component of $\pi$, so we simply refer to it as $\pi$. Further, we note that $u, y$ and $z$ depend upon their initial conditions. We refer to this dependence for the trajectory for $u$ as $u(t,x) = \phi[u_0, z_0](t,x)$ and $z(t) = \psi[u_0,z_0](t).$ We then consider the stochastic process 
\[\pi(u(t),y(t),z(t)) = y(t) + \int_0^{\infty}g(\phi[ u(t), z(t)](s,x_0))ds.\]
Then we use the Ito formula to write down the process for $\pi$: 
\begin{equation}
\begin{aligned}
     d\pi= \,&\nabla \pi(u,y,z)\{du,dy,dz\} +\dfrac{1}{2} \nabla^2\pi(u,y,z)\{du,dy,dz\}
      \\= \,& \nabla \pi(u,y,z)\{f(u,z)dt + \sigma dW,g(u(x_0))dt,h(u(x_0))dt\} \\ \, &+\dfrac{1}{2} \nabla^2\pi(u,y,z)\{f(u,z)dt+\sigma dW,g(u(x_0))dt, h(u(x_0))dt\}.
\end{aligned}
\end{equation}
The process $d\pi$ describes exactly the stochastic phase defined by the isochronal map.

We approximate this stochastic phase, and obtain a closed system by evaluating everything on the stable manifold and replacing $y$ with $\tilde{y}$ as: 
\begin{align*}d\tilde{y} = \nabla \pi(u_*,\tilde{y},z_*)\{f(u_*,z_*)dt+\sigma dW,g(u_*(x_0))dt, h(u_*(x_0))dt\}
\\
+\dfrac{1}{2} \nabla^2\pi(u_*,\tilde{y},z_*)\{f(u_*,z_*)dt+\sigma dW,g(u_*(x_0))dt,h(u_*(x_0))dt\}.
\end{align*}
Note that $g(u_*(x_0)) = h(u_*(x_0)) = f(u_*,z_*)= 0$, since $(u_*,y,z_*)$ is a fixed point of the deterministic system. We are thus left with the SDE
\[d\tilde{y} = \sigma\nabla_u \pi(u_*,\tilde{y},z_*)\{dW\}
+\dfrac{\sigma^2}{2} \nabla_u^2\pi(u_*,\tilde{y},z_*)\{dW\}.
\]
If our noise can be expand as follows 
\[ dW = \sum_{k=1}^{n}v_kdW_k,\] then cross terms vanish since $\mathbb{E}[dW_idW_j]= \delta_{i,j}dt,$ and we obtain 
 \[d\tilde{y} = \sigma \sum_{k=1}^{n}\nabla_u \pi(u_*,\tilde{y},z_*)\{v_k\}dW_k
+\dfrac{\sigma^2}{2} \sum_{k=1}^{n}\nabla_u^2\pi(u_*,\tilde{y},z_*)\{v_k\}dt.
\]
Our task amounts to computing the first and second directional derivatives of $\pi$:
\begin{equation}
    \nabla_u\pi(u_*,\tilde{y},z_*)\{v_k\}= \int_0^{\infty} g'(\phi[u_*,z_*](s,x_0))\nabla_u\phi[u_*,z_*]\{v_k\}(s,x_0) ds.
\end{equation}
\begin{remark}
    Intuitively, we may think that $\nabla_u \phi$ is a matrix multiplying a vector $v$ which yields a vector that we now evaluate at component $x_0$ for each time $s$.
\end{remark}
For the second derivative
\begin{equation}
\begin{aligned}
    \nabla^2_u\pi(u_*,\tilde{y},z_*)\{v_k\}= \int_0^{\infty}(g''(\phi[u_*,z_*](s,x_0)(D_u\phi[u_*,z_*]\{v_k\}(s,x_0))^2 \\
+g'(\phi[u_*,z_*](s,x_0))D^2_u\phi[u_*,z_*]\{v_k\}(s,x_0))ds.
\end{aligned}
\end{equation}
\begin{remark}
    $D_u^2 \phi$ can be thought of as a rank $3$ tensor contracted with the vector $v$ twice yielding another vector evaluated at the component $x_0$ for each time $s$.
\end{remark}
Since $u_*$ and $z_*$ correspond to the stationary solution, these become 
\begin{equation}
\nabla_x\pi(u_*,\tilde{y},z_*)\{v_k\}= \int_0^{\infty}g'(u_*)\nabla_x\phi(s,u_*,z_*)\{v_k\}ds.
\end{equation}
and 
\begin{equation}
  \begin{aligned}
    \nabla^2_x\pi(u_*,\tilde{y})\{v_k\}= \int_0^{\infty}(g''(u_*)(D_u\phi(s,u_*)\{v_k\}[x_0])^2 \\
+g'(u_*)D_u^2\phi(s,u_*)\{v_k\}[x_0])ds.
\end{aligned}
    \label{E:Drift_term}
\end{equation}

Now we compute the directional derivatives of $\phi.$ We do so by considering 
\[D_u\phi(t,u_*,z_*)\{v_k\}=\lim_{\epsilon_0 \to 0} \dfrac{\phi(t,u_* +\epsilon_0v_k,z_*) - \phi(t,u_*,z_*)}{\epsilon_0}.\]
To calculate, we make the ansatz 
\begin{align*} 
\phi(t,u_*+\epsilon_0v_k,z_*) = u_*+\epsilon_0v(t)
\\
z(t)=z_* + z(t)
\end{align*}
where $v(0)=v_k$ and $z(0) = z_*$. Plugging in our ansatz in \eqref{E:Toy_SDE} and Taylor expanding, we find 
\begin{align*} \dfrac{dv(t)}{dt} &= \partial_zf(u_*,z_*)z(t) + D_uf(u_*,z_*)\{v(t)\} +o(1)\\
\dfrac{dz(t)}{dt} &= h'(u_*(x_0))v(x_0,t) +o(1) 
\end{align*} 
Thus the derivative we seek is given by
\begin{equation}
D_u \phi(t,u_*)\{v_k\} = v(t), \begin{aligned}\quad  \dfrac{dv}{dt} &= D_{u,z}f(u_*,z_*)\{(v,z)\}, \quad v(0) = v_k \\
\dfrac{dz}{dt} &= h'(u_*(x_0))v(x_0)  
 , \quad z(0)=0.
\end{aligned}
\label{E:Diffusion_term} \end{equation}
For efficiency, we will refer to $\mathcal{L}(v,z) =D_{u,z}f(u_*,z_*)\{(v,z)\}$. Note the equation is linear so $(-v,-z)$ also solves the equation with initial condition $(-v_k,0)$.

Now we compute the second derivative term
\begin{equation}
    D_u^2 \phi (0,t,u_*)\{v_k\} = \lim_{\epsilon_0 \to \infty} \dfrac{\phi(t,u_*+\epsilon_0v_k) - 2\phi(t,u_*) + \phi(t,u_*-\epsilon_0v_k)}{\epsilon_0^2}.
    \label{E:Second_derivative}
\end{equation}
We plug in the ansatz
\begin{align*}\phi(t,u_*\pm \epsilon_0) =u_*\pm \epsilon_0v(t)+\epsilon_0^2\nu_{\pm}(t) 
\\
z(t) =z_* \pm \epsilon_0 z(t) + \epsilon_0^2 \zeta_{\pm}(t)
\end{align*}
into \eqref{E:Toy_SDE}.
At $O(\epsilon_0)$ 
\begin{align*}\dfrac{dv(t)}{dt} &= \mathcal{L}(v,z) \\ 
\dfrac{dz(t)}{dt} &= h'(u_*(x_0))v(x_0) 
\end{align*}
consistent for what we already had for the first derivative.
At $O(\epsilon_0^2)$ we find that
\begin{equation}
\begin{aligned}
\dfrac{d\nu_{\pm}}{dt} &= 
\mathcal{L}(\nu_{\pm},\zeta_{\pm}) + \dfrac{1}{2}D^2_uf(u_*,z_*)\{v\}+ zD_uD_zf(u_*,z_*)\{v\} +\dfrac{1}{2}\partial_z^2f(u_*,z_*)z^2
\\
\dfrac{d\zeta_{\pm}}{dt} &= h'(u_*(x_0))\nu_{\pm}(x_0)+ h''(u_*(x_0))v(x_0)
.
\end{aligned}
\label{E:PDE_for_second_derivative}
\end{equation}
Note that $(\nu_{+},\zeta_{+})$ and $(\nu_{-},\zeta_{-})$ solve the same equation. We obtain that 
\begin{equation}
D^2_u\phi(t,u_*)\{v_k\} = 2\nu, \begin{aligned}
\quad \dfrac{d\nu}{dt} &= \mathcal{L}(\nu,\zeta) +\dfrac{1}{2} (D_{u,z})^2f(u_*,z_*)\{(v,z)\}, \quad \nu(0) = 0\\
\dfrac{d\zeta}{dt} &= h'(u_*(x_0))\nu(x_0)+ h''(u_*(x_0))v(x_0), \quad \zeta(0) = 0.
\end{aligned}
\label{E:Drift_PDE}
\end{equation}
Note we have already needed to solve for $(v,z)$ in the equation for the first derivative, and the first derivative is used in a forcing term for the equation for the second derivative.

\end{document}